\documentclass[11pt]{article}
\usepackage{tikz}
\usepackage{subfigure}
\usepackage{enumitem}
\usepackage{graphicx,xcolor}
\usepackage{epstopdf}
\usepackage{dsfont}

\usepackage[hyphens]{url}

\usepackage{booktabs} 

\graphicspath{{./figures/}}

\usepackage[onehalfspacing,nodisplayskipstretch]{setspace}

\usepackage{amsfonts} 
\usepackage{bigints}
\usepackage{bm}
\usepackage{amsmath,amssymb,amsbsy,amsthm}

\usepackage{epsfig}
\usepackage{cleveref}
\usepackage{changepage}

\newcommand {\be}[1]{\begin{equation}\label{#1}}
\newcommand {\ee}{\end{equation}}
\newcommand {\bea}{\begin{eqnarray}}
\newcommand {\eea}{\end{eqnarray}}

\newcommand{\mb}[1]{\mbox{\boldmath $#1$}}

\newcommand{\mbs}[1]{{\mbox{\boldmath \scriptsize{$#1$}}}}

\def\texitem#1{\par\smallskip\noindent\hangindent 25pt
               \hbox to 25pt {\hss #1 ~}\ignorespaces}

\newtheorem{theorem}{Theorem}

\newcommand{\appsection}[1]{\let\oldthesection\thesection
  \renewcommand{\thesection}{Electronic Companion \oldthesection}
  \section{#1}\let\thesection\oldthesection}

\theoremstyle{remark}

\newtheorem{example}{Example}

\numberwithin{equation}{section}
\numberwithin{theorem}{section}
\numberwithin{corollary}{section}
\numberwithin{prop}{section}

\graphicspath{{./figs/}}

\definecolor{darkred}{RGB}{139,0,0}
\definecolor{darkgreen}{RGB}{0,139,0}

\begin{document}

\title{Distributionally robust optimization through the lens of submodularity}
\author{Karthik Natarajan\thanks{Engineering Systems and Design, Singapore University of Technology and Design, 8 Somapah Road, Singapore 487372. Email: karthik\_natarajan@sutd.edu.sg} \and Divya Padmanabhan\thanks{School of Mathematics and Computer Science, Indian Institute of Technology, Ponda-403401,
Goa, India. Email: divya@iitgoa.ac.in} \and Arjun Ramachandra\thanks{E207, Decision Sciences Area, Indian Institute of Management Bangalore-560076, India. Email: arjun.ramachandra@iimb.ac.in}
}
\date{May 2025}
\maketitle

\begin{abstract}
Distributionally robust optimization is used to tackle decision making problems under uncertainty where the distribution of the uncertain data is ambiguous. Many ambiguity sets have been proposed for continuous uncertainty that build on convexity and for which the resulting formulations scale polynomially in the number of random variables. However fewer ambiguity sets have been proposed for discrete uncertainty where the exact formulations scale polynomially in the number of random variables. 
Towards this, we define a submodular ambiguity set and showcase its expressive power in modeling both discrete and continuous uncertainty. With discrete uncertainty, we show that a class of  distributionally robust optimization problems is solvable in polynomial time by viewing it through the lens of submodularity. With continuous uncertainty, we show that it is solvable approximately up to an additive error in pseudo-polynomial time. We then focus on a specific class of submodular ambiguity sets where univariate marginal information and bivariate dependence information on the random vector is specified and provide an exact reformulation as a polynomial sized linear program when the uncertainty is discrete and as a polynomial sized semidefinite program when the uncertainty is continuous. We provide numerical evidence of the modeling flexibility and expressive power of the submodular ambiguity set and demonstrate its applicability in two examples: project networks and multi-newsvendor problems. The paper highlights that the submodular ambiguity set is the natural discrete counterpart of the convex ambiguity set and supplements it for continuous uncertainty, both in modeling and computation. \end{abstract}

\maketitle

\section{Introduction}  \label{sec1}
Consider a distributionally robust optimization problem of the form:
\begin{equation} \label{dro}
\begin{array}{rlll}
\displaystyle \inf_{\mbs{x} \in {\cal X}}  \sup_{\mathbb{P} \in {\cal P}}\mathbb{E}_{\mathbb{P}}\left[g(\mb{x},\tilde{\mb{\xi}})
\right],
\end{array}
\end{equation}
where the decision vector $\mb{x}$ is chosen from the set ${\cal X}$ before the true realization of the random vector $\tilde{\mb{\xi}}$ is realized. The probability distribution of $\tilde{\mb{\xi}}$ is denoted by $\mathbb{P}$ and is itself ambiguous. The distribution $\mathbb{P}$ is however assumed to lie in a set of probability distributions denoted by ${\cal P}$, commonly referred to as an \textit{ambiguity set}. The cost incurred for a decision $\mb{x}$ and a realization of the random vector $\tilde{\mb{\xi}} = \mb{\xi}$ is given by $g(\mb{x},\mb{\xi})$. In formulation \eqref{dro}, the decision $\mb{x} \in {\cal X}$ is selected to minimize the worst-case expected cost which is computed over all distributions $\mathbb{P} \in {\cal P}$ and hence is termed a \textit{distributionally robust optimization} problem. When the set ${\cal P} = \{\mathbb{P}\}$ consists of a single distribution, \eqref{dro} reduces to a classical stochastic optimization problem and when ${\cal P} = {\cal P}(\Xi)$ where ${\cal P}(\Xi)$ is the set of all probability distributions with support contained in a closed bounded set ${\Xi}$, \eqref{dro} reduces to a robust optimization problem (see \cite{bentalbook,hertog}). Distributionally robust optimization lies between the two extremes and aims to control for the conservatism of planning for a worst-case realization in robust optimization and planning for a distribution in stochastic optimization that might possibly be mis-specified. Several ambiguity sets have now been proposed in the literature that include the moment ambiguity set \cite{Delage,bkt,wiesemann}, the marginal ambiguity set \cite{chenkarthik,doan1,doan2}, the phi-divergence ambiguity set \cite{hertogbental} and the Wasserstein ambiguity set \cite{kuhnpeyman,murthy,gao}. The recent article of \cite{sor-kuhn} provides an excellent review on the various ambiguity sets and their use in distributionally robust optimization. Broadly, computational tractability is guaranteed when the support $\Xi$ of the uncertainty is a convex set and upper bounds on the expected value of a finite set of convex functions of the random vector are specified in the ambiguity set. Furthermore, such a convex ambiguity set has powerful modeling capabilities. We leave a precise discussion on the details to the next section. When the support $\Xi$ is a discrete set however, the complexity of solving the distributionally robust optimization problem \eqref{dro} typically scales with the size of the set $\Xi$. This in turn might be exponential in the number of random variables. An exception is the marginal ambiguity set where only marginal information on the random vector is specified. The resulting formulations for the marginal ambiguity set for certain cost functions only grows linearly in the number of random variables; see Chapter 2 in \cite{karthikbook}. While the marginal ambiguity set incorporates distributional information on each random variable, it does not explicitly model the dependence among the random variables. A natural question is whether there exist ambiguity sets for discrete uncertainty that allow for modeling both the marginal and dependence information while preserving  computational tractability in solving distributionally optimization problems? Furthermore, can these ambiguity sets be extended to model continuous uncertainty? 

The main contributions of the paper towards answering these questions are described next:

\texitem{(a)} \textit{Modeling uncertainty with a submodular ambiguity set:} We define a submodular ambiguity set where the support of each random variable is specified along with upper bounds (respectively lower bounds) on the expected values of a finite set of submodular (respectively supermodular) functions of the random vector. We illustrate with examples, that marginal and dependence information on the random variables can be modeled using a submodular ambiguity set and compare it to information modeled by existing ambiguity sets in the literature.
\texitem{(b)} \textit{Computational tractability with discrete uncertainty:} When the support of each random variable is discrete and finite we show that the inner supremum in \eqref{dro} is computable in polynomial time under appropriate assumptions on the cost function and the submodular ambiguity set. The result makes use of the ellipsoid method and the polynomial time solvability of submodular function minimization. In turn this helps us identify a class of polynomial time solvable distributionally robust optimization problems by viewing it through the lens of submodularity. These results can be viewed as the discrete counterpart of the known polynomial time complexity result for the convex ambiguity set. We then focus on a specific class of submodular ambiguity sets defined using univariate marginal and bivariate dependence information and develop a polynomial sized linear program when the uncertainty is discrete. 
\texitem{(c)} \textit{Computational tractability with continuous uncertainty:} 
When the support of of each random variable lies in a finite interval, under appropriate assumptions on the cost function and the submodular ambiguity set, we show that the problem is solvable approximately up to additive error in pseudo-polynomial time. This helps us incorporate information on the random variables beyond what is currently modeled by a convex ambiguity set, though at the cost of implementing pseudo-polynomial time algorithms instead of polynomial time algorithms. Interestingly, we show that for a moment based submodular ambiguity set where the means and variances of nonnegative random variables are given along with lower bounds on covariances, one can develop an exact polynomial sized semidefinite program to solve the problem. The corresponding problem when the covariances are exactly specified is known to be NP-hard and thus the approach provides a computationally tractable relaxation.
\texitem{(d)} \textit{Numerical experiments:} We numerically showcase the expressive power and modeling flexibility of the submodular ambiguity set in two applications. In project networks, we use the ambiguity set to incorporate different types of activity duration distributions through discretization along with lower bounds on the correlations among the activities. Using linear programming, we compute the  the worst-case expected project duration in project networks and the corresponding criticality indices and show that the approach provides meaningful estimates. We also consider an extension of the multi-newsvendor problem by incorporating fairness into the classic profit maximization framework. We use the moment based submodular ambiguity set to find the optimal order quantities with semidefinite programming. Numerical results in this application illustrate that incorporating lower bounds on correlations, when available, can provide significant improvement in the performance over distributionally robust optimization models that neglect the information.

The structure of the paper is discussed next. We list the mathematical notations used in the paper in Section \ref{notation}. In Section \ref{sec2}, we review  ambiguity sets that are most related to this paper. In Section \ref{sec3}, we define the submodular ambiguity set and discuss the information that can be modeled using it. We discuss the computational tractability of the submodular ambiguity set with discrete uncertainty in Section \ref{sec4} and with continuous uncertainty in Section \ref{sec5}. Numerical experiments are reported in Section \ref{sec6} before we conclude in Section \ref{sec7}. All proofs are provided in the Appendix A with a review of submodularity provided in Appendix B.

\subsection{Notations}  \label{notation}
We use a nonbold symbol such as $x$ to denote a scalar, bold symbol such as $\mb{x}$ to denote a vector and bold capital symbol such as $\mb{X}$ to denote a matrix. Random numbers and random vectors are denoted with the tilde sign; examples are $\tilde{\xi}$ and $\tilde{\mb{\xi}}$. For a positive integer $N$, we use $[N]$ to denote the set $\{1,2,\ldots,N\}$ and $[N]_{2}$ to denote the set $\{(1,2),(1,3),\ldots,(N-1,N)\}$. Given $p \geq 1$, the $p$-norm of $\mb{x}$ is denoted by $\|\mb{x}\|_p$. The cardinality of a set $\Xi$ is denoted by $|\Xi|$ (possibly infinite) and the closure by $\mbox{cl}(\Xi)$ . Given sets $\Xi_i \subseteq \mathbb{R}$ for $i \in [N]$, the Cartesian product set is given by $\prod_{i \in [N]}\Xi_i = \{({\xi}_1,{\xi}_2,\ldots,{\xi}_N) \ {|} \ {\xi}_1 \in \Xi_1, {\xi}_2 \in \Xi_2, \ldots, {\xi}_N \in \Xi_N\}$. Examples include the Boolean hypercube $\{0,1\}^N$, the hypercube $[0,1]^N$, the real space $\mathbb{R}^N$ and the nonnegative orthant $\mathbb{R}^N_+$. The indicator function of membership in a set $\Xi$ is given by $\mathds{1}_{\mbs{\xi} \in \Xi}$ which takes a value of $1$ if $\mb{\xi}  \in \Xi$ and $0$ if $\mb{\xi}  \not\in \Xi$. The set of real valued symmetric square matrices of size $N \times N$ is denoted by $\mathbb{S}^N$. Given $\mb{X}, \mb{Y} \in \mathbb{S}^N$, we let $\mb{X} \bullet \mb{Y} = \mbox{trace}(\mb{X} \mb{Y})$. We use $\mb{X} \succeq 0$ to denote the matrix $\mb{X}$ is positive semidefinite and $\mb{X} \succeq \mb{Y}$ (respectively $\mb{X} \preceq \mb{Y}$) when $\mb{X}-\mb{Y} \succeq 0$ (respectively $\mb{Y} -\mb{X}\succeq 0$). The vector formed with the diagonal entries of the matrix $\mb{X}$ is denoted by $\mbox{diag}(\mb{X})$. Associated with a random vector $\tilde{\mb{\xi}}$ is a probability distribution $\mathbb{P}$ which we denote by $\tilde{\mb{\xi}} \sim \mathbb{P}$. We use $\mathbb{P}(\cdot)$ to denote the probability of an event, $\mathbb{E}_{\mathbb{P}}[\cdot]$ to denote the expectation with respect to $\mathbb{P}$ and $\text{supp}(\mathbb{P})$ to denote the support of $\mathbb{P}$. We use ${\cal P}(\Xi)$ to denote the set of all probability distributions with support contained in the set $\Xi$. A Dirac measure at $\mb{\xi}$ is denoted by $\delta_{\mbs{\xi}}$. 
The covariance of $\tilde{\xi}$ and $\tilde{\chi}$ is denoted by $\mbox{Cov}(\tilde{\xi},\tilde{\chi})$. The projection of the probability distribution of  $(\tilde{\mb{\xi}},\tilde{\mb{\chi}}) \sim \mathbb{P}$ on $\tilde{\mb{\xi}}$ is given $\mbox{proj}_{\tilde{\mbs{\xi}}}(\mathbb{P})$ and likewise the projection of the ambiguity set of $(\tilde{\mb{\xi}},\tilde{\mb{\chi}}) \sim \mathbb{P} \in {\cal P}$ on  $\tilde{\mb{\xi}}$  is given by  $\mbox{proj}_{\tilde{\mbs{\xi}}}({\cal P})$ .

\section{Ambiguity sets in distributionally robust optimization}  \label{sec2}
Convexity plays a key role in the design of tractable ambiguity sets. We discuss one such ambiguity set that is also key to our work. Given a closed convex set $\Xi \subseteq \mathbb{R}^N$, a set of convex functions $f_l: \Xi \rightarrow \mathbb{R}$ for $l \in [L]$ and a set of scalars $\gamma_l \in \mathbb{R}$ for $l \in [L]$, define the ambiguity set ${\cal P}_{\text{cvx}}$ as follows:
\begin{equation} \label{eq:convex}
\begin{array}{lll}
 {\cal P}_{\text{cvx}}  := \{\mathbb{P} \in {\cal P}(\Xi) \ | \ \mathbb{E}_{\mathbb{P}}[f_l(\tilde{\mb{\xi}})] \leq \gamma_l, \forall l \in [L]\big\},
 \end{array}
 \end{equation}
where all the expectations are assumed to be well-defined with respect to the distributions in the set. We refer to \eqref{eq:convex} as a \textit{convex ambiguity set}. The corresponding distributionally robust optimization problem is given by:
\begin{equation} \label{dro00}
\begin{array}{rlll}
\displaystyle \inf_{\mbs{x} \in {\cal X}} \sup_{\mathbb{P} \in {\cal P}_{\text{cvx}}}\mathbb{E}_{\mathbb{P}}\left[g(\mb{x},\tilde{\mb{\xi}})
\right].
\end{array}
\end{equation}
Under the assumptions: (i) the convex set $\Xi$ has a polynomial time separation oracle, (ii) each of the convex functions $f_l: \Xi \rightarrow \mathbb{R}$ in the ambiguity set has a polynomial time subgradient oracle that returns the function value and its subgradient efficiently and (iii) the cost function is given by $g(\mb{x},{\mb{\xi}}) := \max_{k\in [K]}g_k(\mb{x},\mb{\xi})$ where each function $g_k: {\cal X} \times \Xi \rightarrow \mathbb{R}$ is concave in $\mb{\xi}$ with a polynomial time supergradient oracle that returns the function value and its supergradient with respect to $\mb{\xi}$ efficiently, the inner supremum in \eqref{dro00} for a fixed $\mb{x} \in {\cal X}$ is computable in polynomial time using the ellipsoid method (see theorem 1.5 in \cite{lasserre9}, chapter 3 in \cite{popescuthesis} and proposition 1 in \cite{Delage}). This in turn has important ramifications on the solvability of the distributionally robust optimization problem. When each function $g_k: {\cal X} \times \Xi \rightarrow \mathbb{R}$ is also convex in $\mb{x}$ with a polynomial time subgradient oracle with respect to $\mb{x}$ and the feasible region ${\cal X}$ is a compact convex set with a polynomial time separation oracle, \eqref{dro00} is solvable in polynomial time using the ellipsoid method. In many applications, the functions and sets are structured and conic representable. In these cases, the distributionally robust optimization problem can be reformulated as a conic program that include linear, second order and semidefinite programs as special cases (see \cite{bentalbook,hertog}). The problem is solvable with optimization software such as MOSEK, Gurobi, CPLEX and tailored open source robust optimization packages such as RSOME \cite{zichen}. When some of the decision variables are discrete, the problems are tackled using mixed integer conic programming solvers. Furthermore, for the more challenging class of two stage distributionally robust linear optimization problems with fixed recourse where the number of terms $K$ in the cost function grow exponentially in the size of the problem, one can employ the rich literature on affine and piecewise affine decision rules to obtain tractable conic approximations (see \cite{bentalbook,hertog}).

In addition to the computational tractability, the convex ambiguity set provides great modeling flexibility. We refer the reader to \cite{wiesemann,zichen,meilin} for interesting examples where upper bounds on even degree moments, upper bounds on dispersion measures such as the mean absolute deviation and the directional deviation of linear combinations of random variables are modeled using the convex ambiguity set. In fact, a \textit{cross moment ambiguity set} ${\cal P}_{\text{cm}}$, proposed in \cite{Delage}, and defined as:
\begin{equation}\label{cm}
\begin{array}{lll}
{\cal P}_{\text{cm}} :=  \{\mathbb{P} \in {\cal P}(\Xi) \ | \ \mathbb{E}_{\mathbb{P}}[\tilde{\mb{\xi}}] = \mb{\mu}, \mathbb{E}_{\mathbb{P}}[(\tilde{\mb{\xi}}-\mb{\mu})(\tilde{\mb{\xi}}-\mb{\mu})')] \preceq\mb{\Sigma} \big\},
 \end{array}
 \end{equation}
 where $\Xi$ is a convex set, the mean of the random vector is fixed to $\mb{\mu}$ and the covariance matrix is bounded from above by $\mb{\Sigma}$ in a positive semidefinite order is also modeled using the ambiguity set \eqref{eq:convex} by a lifting argument. Specifically it was shown in \cite{wiesemann} that:
 \begin{equation*} \label{eq:convex1}
\begin{array}{lll}
 {\cal P}_{\text{cm}} = \mbox{proj}_{\tilde{\mbs{\xi}}}(\hat{\cal P}_{\text{cm}})  \mbox{ with }\hat{\cal P}_{\text{cm}} := \{\hat{\mathbb{P}} \in {\cal P}(\tilde{\mb{\xi}} \in \Xi, \tilde{\mb{U}} \succeq (\tilde{\mb{\xi}}-\mb{\mu})(\tilde{\mb{\xi}}-\mb{\mu})')  \ | \ \mathbb{E}_{\hat{\mathbb{P}}}[\tilde{\mb{\xi}}] =\mb{\mu}, \mathbb{E}_{\hat{\mathbb{P}}}[\tilde{\mb{U}}] =\mb{\Sigma}\big\},
 \end{array}
 \end{equation*}
 where $(\tilde{\mb{\xi}},\tilde{\mb{U}}) \sim \hat{\mathbb{P}} \in \hat{\cal P}_{\text{cm}}$ and $\mbox{proj}_{\tilde{\mbs{\xi}}}(\hat{\cal P}_{\text{cm}})$ is the projection of the ambiguity set $\hat{\cal P}_{\text{cm}}$ on the $\tilde{\mb{\xi}}$ sub-vector. The set ${\hat{\cal P}}_{\text{cm}}$ is an example of a convex ambiguity set where the support is given by the intersection of two convex sets since:
   \begin{equation*}
\begin{array}{ll}
{\mb{U}} \succeq ({\mb{\xi}}-\mb{\mu})({\mb{\xi}}-\mb{\mu})' \Leftrightarrow \begin{pmatrix}
    1 & {\mb{\xi}}-\mb{\mu}\\
    {\mb{\xi}}-\mb{\mu} &{\mb{U}}\\
  \end{pmatrix} \succeq 0,
\end{array}
\end{equation*}
and all the expectation constraints are defined on linear functions of the random variables. Incorporating lower bounds on the covariance matrix in a positive semidefinite order is known to be much harder computationally; see remark 1 in \cite{Delage}. Another popular class of ambiguity sets that builds on the convex ambiguity set is the Wasserstein ambiguity set (see \cite{kuhnpeyman,murthy,gao}). These ambiguity sets are characterized by distributions that lie in a Wasserstein ball of a fixed radius around a reference distribution where the reference distribution is typically assumed to be discrete and sample based. When the the support of the uncertainty lies in a convex set, the Wasserstein ambiguity set is computationally tractable.

While continuous uncertainty has been well studied in distributionally robust optimization, in some applications, modeling the data as discrete random variables is natural. Discrete uncertainty is used to represent demand of indivisible goods, failures in power systems, influence relationship among individuals in social networks, labels of images in classification tasks and tokens in large language models. Discretization/binning is also a common technique used to convert any continuous random variable to a discrete random variable. Ambiguity sets that deal with discrete random variables are hence important. One such ambiguity set which allows for the marginal distributions (discrete or continuous) to be either completely or partly specified is the marginal ambiguity set which contains all joint distributions with given marginal information. In the marginal ambiguity set with each random variable taking a finite number of values, the number of possible realizations of the random vector is exponential in the number of random variables. Regardless, under specific assumptions on the cost function and the set ${\cal X}$, \eqref{dro} is solvable in polynomial time using linear optimization; see Chapter 2 in \cite{karthikbook}. However such an ambiguity does not explicitly model the dependence information among the random variables. These ambiguity sets have since been generalized to nonoverlapping and overlapping marginals and in certain cases \eqref{dro} is still efficiently solvable, see Chapter 4 in \cite{karthikbook}. Another tractable ambiguity set that models discrete uncertainty is the phi-divergence ambiguity set (see \cite{hertogbental}). In the input specification of this ambiguity set, all the realizations of the discrete random vector are explicitly listed out and the probabilities of the realizations are allowed to vary. However if the number of realizations is potentially large, these formulations become computationally intractable. The Wasserstein ambiguity set allows for tractable reformulations with discrete uncertainty in specific instances; see the recent work of \cite{selvi}. In this paper, we introduce a class of ambiguity sets that models discrete marginal information and incorporate specific forms of dependence information while preserving computational tractability. Furthermore, we show that these ambiguity sets provide useful modeling capability for continuous uncertainty.

\section{Modeling uncertainty with a submodularity ambiguity set} \label{sec3}
In this section, we discuss modeling of uncertainty, discrete or continuous, with a submodular ambiguity set. Given sets $\Xi_i \subseteq \mathbb{R}$ for each $i \in [N]$, a set of submodular functions $f_l: \prod_{i \in [N]}\Xi_i \rightarrow \mathbb{R}$ for $l \in [L]$ and a set of scalars $\gamma_l \in \mathbb{R}$ for $l \in [L]$, define the ambiguity set ${\cal P}_{\text{sub}}$ as follows:
\begin{equation} \label{eq:submod}
\begin{array}{lll}
 {\cal P}_{\text{sub}}  := \{\mathbb{P} \in {\cal P}(\prod_{i\in [N]}\Xi_i) \ | \ \mathbb{E}_{\mathbb{P}}[f_l(\tilde{\mb{\xi}})] \leq \gamma_l, \forall l \in [L]\big\},
 \end{array}
 \end{equation}
where all the expectations are assumed to be well-defined with respect to the distributions in the set.  We refer to \eqref{eq:submod} as a \textit{submodular ambiguity set}. The submodularity of the functions in the ambiguity set implies they satisfy the condition:
\begin{equation} \label{submodular1}
\begin{array}{rlll}
 f_l(\mb{\xi}) +  f_l(\mb{\chi}) \geq  f_l(\mb{\xi} \wedge \mb{\chi}) + f_l(\mb{\xi} \vee \mb{\chi}), \forall \mb{\xi}, \mb{\chi} \in \prod_{i \in [N]}\Xi_i, \forall l \in [L],
\end{array}
\end{equation}
where $\mb{\xi} \wedge \mb{\chi} = (\min(\xi_1,\chi_1),\ldots,\min(\xi_N,\chi_N))$ and $\mb{\xi} \vee \mb{\chi} = (\max(\xi_1,\chi_1),\ldots,\max(\xi_N,\chi_N))$. When the functions are twice differentiable, this condition is equivalent to:
\begin{equation} \label{submodular1a}
\begin{array}{rlll}
  {\displaystyle \frac{\partial^2}{\partial \xi_i \partial \xi_j} f_l(\mb{\xi}) \leq 0}, \forall \mb{\xi} \in  \prod_{i \in [N]}\Xi_i, \forall i \neq j, \forall l \in [L].
  \end{array}
\end{equation}
In Appendix B, we review submodularity, submodular function minimization and comonotonic random variables which will be crucial for this paper. Much of this material is standard and the reader can refer to excellent articles and books on this topic; see \cite{Lovasz,topkis,bach}. The similarity to the convex ambiguity set is immediate where convex functions are replaced with submodular functions and the convex set $\Xi \subseteq \mathbb{R}^N$ is replaced with the Cartesian product set $\prod_{i\in [N]}\Xi_i \subseteq \mathbb{R}^N$. While the ambiguity set in \eqref{eq:submod} specifies upper bounds on the expected value of submodular functions, equivalently, we can specify lower bounds on the expected value of supermodular functions in ${\cal P}_{\text{sub}}$. This follows from the fact that a function $f_l: \prod_{i \in [N]}\Xi_i \rightarrow \mathbb{R}$ is supermodular if and only if $-f_l: \prod_{i \in [N]}\Xi_i \rightarrow \mathbb{R}$ is submodular. The sign of the inequalities in \eqref{submodular1} and \eqref{submodular1a} are thus reversed in the definition of supermodular functions. We next illustrate some types of information on the uncertainty that can be modeled with a submodular ambiguity set.

\subsection{Univariate information}
 Give a finite set of arbitrary univariate functions $f_{i,l}: \Xi_i \rightarrow \mathbb{R}$ and scalars $\underline{\gamma}_{i,l}, \overline{\gamma}_{i,l} \in \mathbb{R}$ that satisfy $\underline{\gamma}_{i,l} \leq \overline{\gamma}_{i,l}$ for each $l \in [L_i]$, $i \in [N]$, define the \textit{univariate ambiguity set} ${\cal P}_{\text{uni}}$ as follows:
\begin{equation} \label{1}
\begin{array}{lll}
 {\cal P}_{\text{uni}} := \{\mathbb{P} \in {\cal P}(\prod_{i \in [N]} \Xi_i) \ | \ \mathbb{E}_{\mathbb{P}}[f_{i,l}(\tilde{\xi}_i)] \in [\underline{\gamma}_{i,l},\overline{\gamma}_{i,l}], \forall l \in [L_i], \forall i \in [N]\big\}.
 \end{array}
 \end{equation}
 Since all univariate functions are submodular by definition, by rewriting the constraints as $\mathbb{E}_{\mathbb{P}}[f_{i,l}(\tilde{\xi}_i)] \leq \overline{\gamma}_{i,l}$ and  $\mathbb{E}_{\mathbb{P}}[-f_{i,l}(\tilde{\xi}_i)]\leq -\underline{\gamma}_{i,l}$, we see that ${\cal P}_{\text{uni}}$ is a submodular ambiguity set. By setting $\underline{\gamma}_{i,l} = -\infty$ (likewise $\overline{\gamma}_{i,l} = +\infty$), we allow for specifying only an upper bound (likewise lower bound) on $\mathbb{E}_{\mathbb{P}}[f_{i,l}(\tilde{\xi}_i)]$ and by setting $\underline{\gamma}_{i,l} = \overline{\gamma}_{i,l}$, we allow for precisely fixing the value in the ambiguity set. The set ${\cal P}_{\text{uni}}$ allows to exactly specify or bound terms such as the marginal probability ${\mathbb{P}}(\tilde{\xi}_i \leq t)$ for a given $t$, marginal moment $\mathbb{E}_{\mathbb{P}}[\tilde{\xi}_i^l]$ for a given $l$, marginal dispersion measure $\mathbb{E}_{\mathbb{P}}[|\tilde{\xi}_i-m_{i,l}|^l]$ for a given $l$ and $m_{i,l}$ or the moment generating function $\mathbb{E}_{\mathbb{P}}[e^{t\tilde{\xi}_i}]$ for a given $t$. While ${\cal P}_{\text{uni}}$ allows for multiple such terms, only finitely many terms are allowed. When the support of each random variable is discrete and finite, ${\cal P}_{\text{uni}}$ can incorporate the full marginal distribution and this reduces to the Fr\'{e}chet set of joint discrete distributions with finite marginal support. Solving the inner problem in  \eqref{dro} with this ambiguity set then reduces to solving an instance of a multimarginal optimal transport problem (see \cite{pass,alt}). A closely related ambiguity set is given by:
     \begin{equation}
\begin{array}{lll}
 {\cal Q}_{\text{uni}} :=\{\mathbb{P} \in {\cal P}(\prod_{i \in [N]} \Xi_i) \ | \ \mathbb{E}[\sum_{i \in [N]}f_{i,l}(\tilde{\xi}_i)] \in [\underline{\gamma}_{l},\overline{\gamma}_{l}], \forall l \in [L]\big\},
 \end{array}
 \end{equation}
 where $f_{i,l}: \Xi_i \rightarrow \mathbb{R}$ are arbitrary univariate functions for all $l$, $i$ and the scalars satisfy $\underline{\gamma}_{l} \leq \overline{\gamma}_{l}$ for all $l$. Such an ambiguity set incorporates information on terms such as $\mathbb{E}_{\mathbb{P}}[\|\tilde{\mb{\xi}}\|_p^p]$ or $\mathbb{E}_{\mathbb{P}}[\|\tilde{\mb{\xi}}-\mb{m}\|_p^p]$.
\subsection{Bivariate information}
 Give a finite set of bivariate submodular functions $f_{i,j,l}: \Xi_i \times \Xi_j \rightarrow \mathbb{R}$ and scalars $\gamma_{i,j,l}\in \mathbb{R}$ for $l \in [L_{i,j}]$, $(i,j) \in [N]_2$, define the \textit{bivariate submodular ambiguity set} ${\cal P}_{\text{bi-sub}}$ as follows:
 \begin{equation}
\begin{array}{lll} \label{2a}
 {\cal P}_{\text{bi-sub}} := \{\mathbb{P} \in {\cal P}(\prod_{i \in [N]} \Xi_i) \ | \ \mathbb{E}[f_{i,j,l}(\tilde{\xi}_i,\tilde{\xi}_j)] \leq \gamma_{i,j,l}, \forall l \in [L_{i,j}], \forall (i,j) \in [N]_2\big\}.
 \end{array}
 \end{equation}
Unlike the univariate ambiguity set, the inequality constraint matters in a bivariate submodular ambiguity set since not all bivariate functions are submodular. Intersecting ${\cal P}_{\text{uni}}$ or ${\cal Q}_{\text{uni}}$ with ${\cal P}_{\text{bi-sub}}$ also gives a submodular ambiguity set. Marginal information on the random variables such as the shape, symmetry, tail behavior, skewness or kurtosis can be modeled using ${\cal P}_{\text{uni}}$ while dependence information among pairs of random variables can be modeled using ${\cal P}_{\text{bi-sub}}$. We provide three such examples next. 
\begin{example}[Bernoulli with lower bounds on covariances]
Given a Bernoulli random vector with known marginal probabilities and lower bounds on the pairwise probabilities of the random variables taking a value 1, consider the ambiguity set:
\begin{equation}
\begin{array}{lll}\label{bergood}
 \{\mathbb{P} \in {\cal P}(\{0,1\}^N) \ | \ \mathbb{P}(\tilde{\xi}_i= 1) = {p}_i, \forall i \in [N],  \mathbb{P}(\tilde{\xi}_i= 1,\tilde{\xi}_j= 1) \geq p_{i,j}, \forall (i,j) \in [N]_2\big\}.
 \end{array}
 \end{equation}
 This ambiguity set models lower bounds on the covariances term by term since $\mbox{Cov}(\tilde{\xi}_i,\tilde{\xi}_j) = \mathbb{P}(\tilde{\xi}_i= 1,\tilde{\xi}_j= 1)-\mathbb{P}(\tilde{\xi}_i= 1)\mathbb{P}(\tilde{\xi}_j= 1)\geq p_{i,j}-p_ip_j $. The set is defined using univariate functions of the form $\mathds{1}_{{\xi}_i = 1}$ and bivariate supermodular functions of the form $\mathds{1}_{{\xi}_i \geq 1, {\xi}_j \geq 1}$ where lower bounds on the expected value of the supermodular functions are given. Hence the ambiguity set is a special case of ${\cal P}_{\text{uni}} \cap {\cal P}_{\text{bi-sub}}$. When $p_i \in [0,1]$ for all $i \in [N]$ and $p_{i,j} = \max(0,p_i+p_j-1)$, the ambiguity set contains all dependent Bernoulli random vectors with fixed marginal probabilities. When $p_{i,j} = \min(p_i,p_j)$, the ambiguity set allows for only the comonotonic (perfectly positively dependent) Bernoulli random vector with fixed marginal probabilities; see Appendix B for a discussion on comonotonic random vectors and the connection to submodularity. As the input parameter $p_{i,j}$ is varied, the ambiguity set controls for the minimum allowable covariance between the random variables $\tilde{\xi}_i$ and $\tilde{\xi}_j$ in the ambiguity set. This ambiguity set has been considered in \cite{BorosScozzari2014} where the authors find polynomial time computable bounds on the probability of the union of random events. A closely related ambiguity set that is not submodular is given by:
 \begin{equation}
\begin{array}{lll} \label{berbad}
 \{\mathbb{P} \in {\cal P}(\{0,1\}^N) \ | \ \mathbb{P}(\tilde{\xi}_i= 1) = {p}_i, \forall i \in [N],  \mathbb{P}(\tilde{\xi}_i= 1,\tilde{\xi}_j= 1) = p_{i,j}, \forall (i,j) \in [N]_2\big\},
 \end{array}
 \end{equation}
 where the pairwise probabilities of the random variables jointly taking a value 1 is precisely specified. Unfortunately testing feasibility or optimizing over the set of distributions \eqref{berbad} is known to be NP-complete; see \cite{pitowsky}. The submodular ambiguity set \eqref{bergood} provides a tractable relaxation of  \eqref{berbad}.

 \end{example}
 \begin{example}[Lower bounds on bivariate tail probabilities or cross moments]
Given a random vector with known univariate marginal information and lower bounds on bivariate tail probabilities, consider the ambiguity set:
\begin{equation} \label{orth}
\begin{array}{lll}
 \{\mathbb{P} \in {\cal P}(\prod_{i \in [N]} \Xi_i) \ | \ \mathbb{E}_{\mathbb{P}}[f_{i,l}(\tilde{\xi}_i)] \in [\underline{\gamma}_{i,l},\overline{\gamma}_{i,l}], \forall l \in [L_i], \forall i \in [N], \\
\quad \quad \quad \quad \quad \quad \quad \quad  \mathbb{P}(\tilde{\xi}_i \geq \xi_i, \tilde{\xi}_j \geq \xi_{j}) \geq \gamma_{i,j}(\xi_{i},\xi_{j}), \forall (\xi_i, \xi_j) \in {\Xi}_{i,j},  \forall (i,j) \in [N]_2\big\},
 \end{array}
 \end{equation}
 where ${\Xi}_{i,j}$ is a finite subset of $\Xi_i \times \Xi_j$ for each $(i,j) \in [N]_2$. Here we use bivariate supermodular functions of the form $\mathds{1}_{{\xi}_i \geq s, {\xi}_j \geq t}$ where lower bounds on the expected value of the supermodular functions are given. When each random variable is assumed to be discrete with finite support and ${\Xi}_{i,j} = \Xi_i \times \Xi_j$, the ambiguity set models the univariate marginal distributions along with bivariate stochastic orders such as the bivariate orthant dependence order and the bivariate concordance order (see \cite{mullerbook}). Such stochastic orders are known to be stronger in capturing dependence relationships among random variables than covariances. A related submodular ambiguity set that models lower bounds on the cross moments is given by:
  \begin{equation} \label{cov}
\begin{array}{lll}
 \{\mathbb{P} \in {\cal P}(\prod_{i \in [N]} \Xi_i) \ | \ \mathbb{E}_{\mathbb{P}}[f_{i,l}(\tilde{\xi}_i)] \in [\underline{\gamma}_{i,l},\overline{\gamma}_{i,l}], \forall l \in [L_i], \forall i \in [N], \mathbb{E}[\tilde{\xi}_i \tilde{\xi}_j] \geq \gamma_{i,j}, \forall (i,j) \in [N]_2\big\}.
 \end{array}
 \end{equation}
 Here we use bivariate supermodular functions of the form ${\xi}_i {\xi}_j$ instead of $\mathds{1}_{{\xi}_i \geq s, {\xi}_j \geq t}$. While \eqref{orth} models stronger bivariate dependence relationships than \eqref{cov}, it comes at the cost of potentially more information having to be specified in the ambiguity set.
 \end{example}
 \begin{example}[Upper bounds on pairwise Wasserstein distances]
 Given a random vector with known univariate marginal information and upper bounds on the moments of the absolute difference of pairs of random variables, consider the ambiguity set:
  \begin{equation} \label{eq:abs}
\begin{array}{lll}
 \{\mathbb{P} \in {\cal P}(\prod_{i \in [N]} \Xi_i) \ | \ \mathbb{E}_{\mathbb{P}}[f_{i,l}(\tilde{\xi}_i)] \in [\underline{\gamma}_{i,l},\overline{\gamma}_{i,l}], \forall l \in [L_i], \forall i \in [N], \\
 \quad \quad \quad \quad \quad \quad \quad \quad \mathbb{E}_{\mathbb{P}}[|\tilde{\xi}_i-\tilde{\xi}_j|^{p_l}] \leq \gamma_{i,j,l}, \forall l \in [L_{i,j}],\forall (i,j) \in [N]_2 \big\}.
 \end{array}
 \end{equation}
Given scalars $p_l \geq 1$, this is an example of a submodular ambiguity set, since each bivariate function $|{\xi}_i-{\xi}_j|^{p_l}$ is convex and hence submodular; see Appendix B. The ambiguity set bounds the Wasserstein distances between every pair of random variables while incorporating univariate marginal information. If we assume the support of the random vector lies in a convex set and each univariate function $f_{i,l}$ is convex with only upper bounds specified,  this is also a convex ambiguity set. In the submodular ambiguity set \eqref{eq:abs} however, we do not put restrictions on the  marginal support nor assume convexity of the univariate functions.
\end{example}
\subsection{Multivariate information}
Here, we discuss submodular ambiguity sets where multivariate information is explicitly modeled.

\begin{example} [Upper bounds on expected norm of a nonnegative random vector] Given $\prod_{i \in [N]} \Xi_i\subseteq \mathbb{R}^N_+$ and upper bounds $\gamma_l$ on the expected $p_l$-norms of the random vector, consider the ambiguity set:
   \begin{equation} \label{3}
\begin{array}{lll}
 \{\mathbb{P} \in {\cal P}(\prod_{i \in [N]} \Xi_i\subseteq \mathbb{R}^N_+) \ | \ \mathbb{E}[||\tilde{\mb{\xi}} ||_{p_l}] \leq \gamma_{l}, \forall l \in [L]\big\}.
 \end{array}
 \end{equation}
 When $p_l \geq 1$, the norm function $||{\mb{\xi}} ||_{p_l}$ is submodular over the nonnegative orthant; see \cite{singla}. Hence this is an example of a submodular ambiguity set. When the support of the random vector is a convex set, this is also a convex ambiguity set.  It should be noted however that an upper bound constraint on the expected norm of a translated random vector of the form $\mathbb{E}[||\tilde{\mb{\xi}} - \mb{m}_l||_{p_l}] \leq \gamma_l$ can be modeled using a convex ambiguity set but not with a submodular ambiguity set since the function $||{\mb{\xi}} - \mb{m}_l||_{p_l}$ is not submodular.
 \end{example}

 \begin{example} [Upper bounds on expected value of concave functions of sum of random variables] Given univariate concave functions $f_l: \mathbb{R} \rightarrow \mathbb{R}$ for $l \in [L]$ and scalars $\gamma_l$ for $l \in [L]$, consider an ambiguity set given by:
  \begin{equation} \label{4}
\begin{array}{lll}
 \{\mathbb{P} \in {\cal P}(\prod_{i \in [N]} \Xi_i) \ | \ \mathbb{E}_{\mathbb{P}}\left[f_l\left(\sum_{i \in [N]} \tilde{\xi}_i\right)\right] \leq \gamma_l, \forall l\in [L] \big\}.
\end{array}
\end{equation}
This is an example of a submodular ambiguity set. Such functions include $\min_{r \in [R]}(a_r\sum_{i \in [N]} {\xi}_i+b_r)$ and $\ln(\sum_{i \in [N]} {\xi}_i)$ defined on the appropriate domain.
\end{example}

\begin{example} [Lower bounds on orthant probabilities] Given a set of vectors ${\mb{\xi}}_{l} \in \prod_{i \in [N]} \Xi_i$  and scalars $\gamma_l \in [0,1]$ for $l \in [L]$, consider the ambiguity set:
  \begin{equation} \label{5}
\begin{array}{lll}
 \{\mathbb{P} \in {\cal P}(\prod_{i \in [N]} \Xi_i) \ | \ \mathbb{P}(\tilde{\mb{\xi}} \geq {\mb{\xi}}_{l}) \geq \gamma_l, \forall l \in [L]\big\}.
 \end{array}
 \end{equation}
 This ambiguity set gives lower bounds on the orthant (tail) probabilities and is an example of a submodular ambiguity set using supermodular functions of the form $\mathds{1}_{{\mbs{\xi}} \geq {\mbs{\xi}}_{l}}$.
  \end{example}
While submodular and supermodular functions have been previously used in robust and distributionally robust optimization problems; see \cite{iancu,agrawal,chenkarthik,jin,mak}, to the best of our knowledge, the use of such functions in explicitly constructing ambiguity sets has not been explored thus far. In the following sections, we discuss the computational tractability of the submodular ambiguity set and its application to distributionally robust optimization. 
\section{Discrete uncertainty}  \label{sec4}
In this section, we provide computationally tractable models for the distributionally robust optimization problem with a submodular ambiguity set when the uncertainty is discrete. Specifically, we consider the problem:
\begin{equation} \label{dro10}
\begin{array}{rlll}
\displaystyle \inf_{\mbs{x} \in {\cal X}} \sup_{\mathbb{P} \in {\cal P}_{\text{sub}}}\mathbb{E}_{\mathbb{P}}\left[g(\mb{x},\tilde{\mb{\xi}})
\right],
\end{array}
\end{equation}
and show that it is solvable in polynomial time under fairly general assumptions on the cost function and the submodular ambiguity set.
\subsection{Polynomial time solvability using the ellipsoid method} \label{sec4a}
Towards solving \eqref{dro10}, we first focus on the problem of computing the inner supremum. Consider the problem:
\begin{equation} \label{droiner}
\begin{array}{rlll}
 \rho^* := {\displaystyle\sup_{\mathbb{P} \in {\cal P}_{\text{sub}}}\mathbb{E}_{\mathbb{P}}\left[g(\tilde{\mb{\xi}})
\right]} = \sup \left\{ \mathbb{E}_{\mathbb{P}}[g(\tilde{\mb{\xi}})] \ | \ \mathbb{P} \in {\cal P}(\prod_{i\in [N]}\Xi_i),  \mathbb{E}_{\mathbb{P}}[f_l(\tilde{\mb{\xi}})] \leq \gamma_l, \forall l \in [L]\big\}
\right\},
\end{array}
\end{equation}
where $g: \prod_{i \in [N]} \Xi_i \rightarrow \mathbb{R}$ is the cost function and $\rho^*$ is the optimal value.
We make the following assumptions on the ambiguity set and the cost function: \\ 
\texttt{Assumption (A1):} The set $\Xi_i \subset \mathbb{R}$ is a discrete, finite set for each $i \in [N]$.\\
\texttt{Assumption (A2):} For each $l \in [L]$, $f_l: \prod_{i \in [N]} \Xi_i \rightarrow \mathbb{R}$ is a submodular function with a polynomial time evaluation oracle. \\
\texttt{Assumption (A3):} The ambiguity set ${\cal P}_{\text{sub}}$ is nonempty. \\
\texttt{Assumption (A4):} The cost function is given by $g({\mb{\xi}}) := \max_{k \in [K]}g_k(\mb{\xi})$ where for each $k \in [K]$, the function $g_k: \prod_{i \in [N]} \Xi_i \rightarrow \mathbb{R}$ is a supermodular function with a polynomial time evaluation oracle.\\
This brings us to the first theorem of the paper.
\begin{theorem} \label{thm:submod}
Suppose assumptions (A1)-(A4) hold. Then $\rho^*$ is computable in polynomial time.  
\end{theorem}
\begin{proof}
See Appendix.
\end{proof}

We make a few remarks about Theorem \ref{thm:submod} and its proof next.
\texitem{(a)} The optimal value $\rho^*$ is computable in time polynomial in $N$ (number of random variables), $L$ (number of submodular functions defining the ambiguity set), $K$ (number of piecewise terms in the objective function), $\max_{i\in[N]}|\Xi_i|$ (the maximum cardinality of any marginal support) and the function evaluation time of the oracles.
It makes use of linear programming duality, the ellipsoid method and the polynomial time solvability of submodular function minimization. The result closely mirrors that for the convex ambiguity set discussed in Section \ref{sec2} where the objective function is given by the maximum of concave functions of the random vector instead (see theorem 1.5 in \cite{lasserre9}, chapter 3 in \cite{popescuthesis} and proposition 1 in \cite{Delage}). We believe the value of Theorem \ref{thm:submod} lies in identifying new polynomial time solvable instances of  (\ref{droiner}). As an illustration, consider the ambiguity set previously discussed in Example 3 which prescribes upper bounds on the Wasserstein distances between pairs of random variables and incorporates univariate marginal information on the discrete random variables. Theorem \ref{thm:submod} illustrates that $\rho^*$ is polynomial time computable for this ambiguity set. To the best of our knowledge, this is not covered by current results in distributionally robust optimization.  
\texitem{(b)} The proof of Theorem \ref{thm:submod} makes use of the ellipsoid method which helps with the theoretical analysis. However, the ellipsoid method is not easy to implement in practice. As we will show in the next section, for structured submodular ambiguity sets that also includes Example 3, it is possible to obtain compact linear programs which makes it practical to use. Similar developments have occurred over the past decade for the convex ambiguity set where conic reformulations have been developed for structured convex ambiguity sets. 
\texitem{(c)} The proof of Theorem \ref{thm:submod} makes use of duality. A natural question is to characterize  the extremal primal distribution which attains the bound $\rho^*$. When the submodular ambiguity set incorporates the full marginal distributions with $\tilde{\xi}_i \sim \mathbb{P}_i$ for each $i \in [N]$, a feasible distribution in the set is given by the comonotonic random vector with perfect positive dependence where $\tilde{\mb{\xi}}^c \sim\mathbb{P}^c(\mathbb{P}_1,\ldots,\mathbb{P}_N)$; see Appendix B. A well known extremal characterization of the comonotonic random vector (see \cite{tchen,ruschbook}) is given by:
\begin{equation*} \label{droinfer}
\begin{array}{rlll}
\displaystyle \sup_{\mathbb{P} \in {\cal P}_F(\mathbb{P}_1,\ldots,\mathbb{P}_N)}\mathbb{E}_{\mathbb{P}}\left[g(\tilde{\mb{\xi}})\right] = \mathbb{E}_{\mathbb{P}^c(\mathbb{P}_1,\ldots,\mathbb{P}_N)}\left[g(\tilde{\mb{\xi}}^c)\right], \forall \mbox{supermodular } g \mbox{ such that expectations exists},
\end{array}
\end{equation*}
where ${\cal P}_F(\mathbb{P}_1,\ldots,\mathbb{P}_N)$ is the Fr\'{e}chet set of joint distributions with fixed marginals. In our context, the objective function is given by the maximum of supermodular functions which is not supermodular in general and so the extremal distribution does not have to be comonotonic. The next example illustrates that the dependence structure in the extremal distribution is related to the specific form of the objective function and the submodular ambiguity set.

\begin{example}[Extremal distribution for maximum of two Bernoulli random variables] Let $N = 2$ and $g({{\xi}_1},{{\xi}_2}) = \max({\xi}_1,{\xi}_2)$. Given the Bernoulli ambiguity set in \eqref{bergood}, consider the problem of computing the maximum expected value:
\begin{equation*} \label{droinerex}
\begin{array}{rlll}
\displaystyle \max \left\{\mathbb{E}_{\mathbb{P}}\left[\max(\tilde{\xi}_1,\tilde{\xi}_2)\right] \ \big{|} \  \mathbb{P} \in {\cal P}(\{0,1\}^N), \mathbb{P}(\tilde{\xi}_1= 1) = {p}_1, \mathbb{P}(\tilde{\xi}_2= 1) = {p}_2, \mathbb{P}(\tilde{\xi}_1= 1,\tilde{\xi}_2= 1) \geq p_{1,2}\right\}.
\end{array}
\end{equation*}
This ambiguity set is nonempty when $p_1,p_2 \in [0,1]$ and $p_{1,2} \leq \min(p_1,p_2)$. We can reformulate the problem by using a single decision variable for the probability $p = \mathbb{P}(\tilde{\xi}_1=1,\tilde{\xi}_2= 1)$ as follows:
\begin{equation*} \label{droinerex}
\begin{array}{rlll}
\displaystyle \max \left\{p_1+p_2-p \ \big{|} \ p_1-p \geq 0, p_2-p \geq 0, 1-p_1-p_2+p \geq 0, p \geq 0, p \geq p_{1,2}\right\},
\end{array}
\end{equation*}
where $\mathbb{P}(\tilde{\xi}_1=1,\tilde{\xi}_2= 0) = p_1-p$, $\mathbb{P}(\tilde{\xi}_1=0,\tilde{\xi}_2= 1) = p_2-p$ and $\mathbb{P}(\tilde{\xi}_1=0,\tilde{\xi}_2= 0) = 1-p_1-p_2+p$. The objective function is given by $\mathbb{E}_{\mathbb{P}}[\max(\tilde{\xi}_1,\tilde{\xi}_2)] = p_1+p_2-p$. The optimal value of $p$ is given by $p^* = \max(0,p_1+p_2-1,p_{12})$. This implies:
\begin{enumerate}
    \item[]{(i)} When $p_{1,2} \leq \max(0,p_1+p_2-1)$, we have $p^* = \max(0,p_1+p_2-1) $ and the extremal distribution is given by the countermonotonic distribution with perfectly negatively dependent random variables.
    \item[]{(ii)} When $p_{1,2} = p_1p_2$, we have $p^* = p_1p_2$ and the extremal distribution is given by the independent distribution. 
    \item[]{(iii)}
    When $p_{12} = \min(p_1,p_2)$, , we have $p^* = \min(p_1,p_2)$ and the extremal distribution is the comonotonic distribution with perfectly positively dependent random variables.
\end{enumerate} Thus for every $p_{1,2} \in [\max(0,p_1+p_2-1),\min(p_1,p_2)]$, the inequality constraint in the ambiguity set is tight at optimality for the max objective function. The extremal distribution hence ranges from perfectly negatively dependent to independent to perfectly positively dependent Bernoulli random variables based on the specified parameter $p_{1,2}$. On the other hand, if we solve the problem over the submodular ambiguity set: \begin{equation*} \label{droinerex}
\begin{array}{rlll}
\displaystyle \max \left\{\mathbb{E}_{\mathbb{P}}\left[\max(\tilde{\xi}_1,\tilde{\xi}_2)\right] \ \big{|} \  \mathbb{P} \in {\cal P}(\{0,1\}^N), \mathbb{P}(\tilde{\xi}_1= 1) = {p}_1, \mathbb{P}(\tilde{\xi}_2= 1) = {p}_2\right\},
\end{array}
\end{equation*}
the extremal distribution will always be the countermonotonic distribution. The example thus illustrates the value of modeling dependence information, when available, using the submodular ambiguity set. 
\end{example}
In the previous example, the objective function is itself submodular and the inequality constraint in the ambiguity set is always tight at the extremal distribution for values of $p_{1,2} \in [\max(0,p_1+p_2-1),\min(p_1,p_2)]$. One might conjecture that the inequality constraints are tight for all submodular objective functions. However, the next example shows that this is not always true.

\begin{example}[Inequality constraints in the ambiguity set are not always tight for submodular objective] We consider the Bernoulli ambiguity sets in \eqref{bergood} and \eqref{berbad} with $p_{i,j}=p_ip_j$ where each $p_i \in [0,1]$. Both the ambiguity sets are nonempty since \eqref{berbad} is the set of pairwise independent Bernoulli random vectors with fixed marginals. Let  $g({\mb{\xi}}) = \max_{(i,j)\in [N]_2} (\xi_i+\xi_j)$ which is a submodular function. Set $N = 5$ and $p_1 = p_2 = \alpha$ and $p_3 = p_4 = p_5 = \beta$ and compute the bounds by varying $\alpha$ and $\beta$. Let $\rho_{\text{eq}}^*$ be the maximum expected value obtained using the ambiguity set \eqref{berbad} (equality constraints) and $\rho_{\text{rel}}^*$ be the maximum expected value obtained using the submodular ambiguity set \eqref{bergood} (relaxed inequality constraints). Clearly $\rho_{\text{rel}}^* \geq \rho_{\text{eq}}^*$. In Figure \ref{fig:alphabetacomb}, we plot the percentage difference for several values of $\alpha$ and $\beta$, clearly highlighting that the bounds can be strictly different with $\rho_{\text{rel}}^* > \rho_{\text{eq}}^*$. 

\begin{figure}[hbtp!]
\begin{center}
 \includegraphics[scale=0.625]{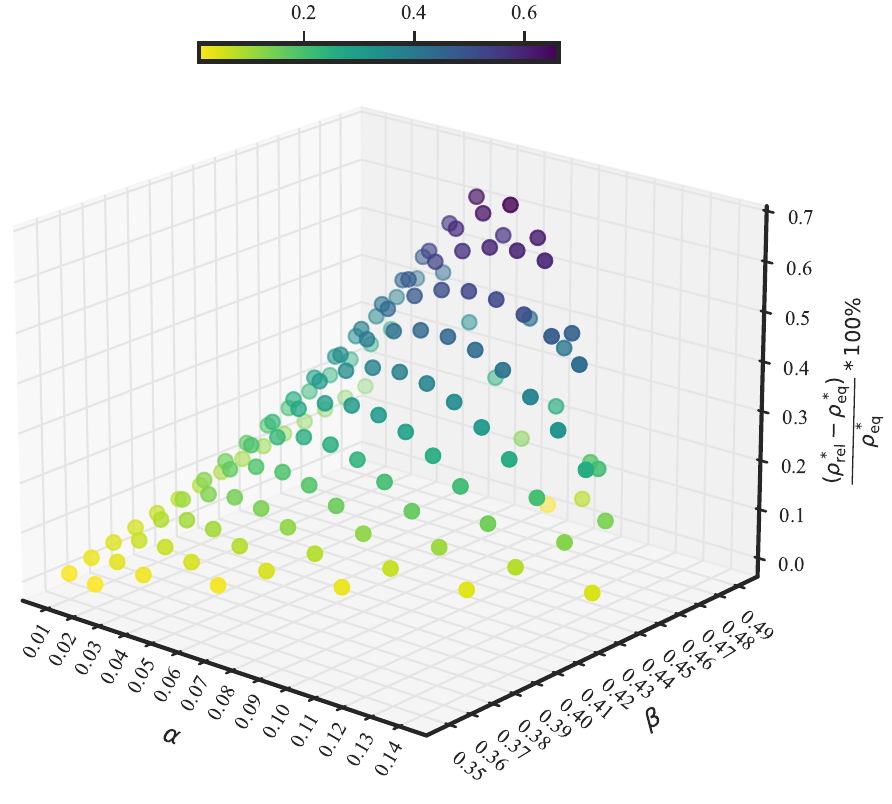}
 \caption{Percentage relative difference in objective values given by $\frac{(\rho_{\text{rel}}^*- \rho_{\text{eq}}^*)}{\rho_{\text{eq}}^*}*100\%$ for varying values of $\alpha,\beta$.}
 \label{fig:alphabetacomb}
 \end{center}
\end{figure}
\end{example}
The next theorem builds on Theorem \ref{thm:submod} to identify conditions that guarantee polynomial time solvability of the distributionally robust optimization problem.

\begin{theorem} \label{thm:submod1}
Consider the distributionally robust optimization problem in \eqref{dro10}. Suppose assumptions (A1)-(A3) hold. Let $g(\mb{x},{\mb{\xi}}) := \max_{k\in [K]}g_k(\mb{x},\mb{\xi})$ where each function $g_k: {\cal X} \times \prod_{i \in [N]}\Xi_i \rightarrow \mathbb{R}$ is convex in $\mb{x}$ with a polynomial time subgradient oracle with respect to $\mb{x}$ and supermodular in $\mb{\xi}$ with a polynomial time evaluation oracle. When the feasible region ${\cal X}$ is a compact convex set with a polynomial time separation oracle, the  distributionally robust optimization problem \eqref{dro10} is solvable in polynomial time.
\end{theorem}
\begin{proof}
See Appendix.
\end{proof}

A similar polynomial time complexity result exists for the convex ambiguity set that makes use of the  ellipsoid method as discussed in Section \ref{sec2}. Theorem \ref{thm:submod1} extends this to discrete uncertainty. 

\subsection{Polynomial sized linear program with univariate and bivariate information}
 \label{sec4b}
In this section, we develop a linear program for a class of submodular ambiguity sets where univariate and bivariate marginal information on the discrete random vector is specified. Here, we focus on the problem of computing the tightest upper bound on the expected value of a piecewise quadratic function of the random vector given by:
\begin{equation} \label{droiner00}
\begin{array}{rlll}
\displaystyle \rho^* := \sup_{\mathbb{P} \in {\cal P}_{\text{uni}} \cap {\cal P}_\text{bi-sub}}\mathbb{E}_{\mathbb{P}}\left[\max_{k \in [K]}\left(\tilde{\mb{\xi}}'\mb{A}_k\tilde{\mb{\xi}} + \mb{b}_k'\tilde{\mb{\xi}} + c_{k}\right)\right].
\end{array}
\end{equation}
The objective function in \eqref{droiner00} is piecewise quadratic with given symmetric matrices $\{\mb{A}_1,\ldots,\mb{A}_K\} \subseteq \mathbb{S}^N$, vectors $\{\mb{b}_1,\ldots,\mb{b}_K\} \subseteq \mathbb{R}^N$ and scalars $\{c_1,\ldots,c_K\} \subseteq \mathbb{R}$. The ambiguity set is given by:
\begin{equation} \label{droinera}
\begin{array}{rlll}
{\cal P}_{\text{uni}} \cap {\cal P}_{\text{bi-sub}} = \{\mathbb{P} \in {\cal P}(\prod_{i \in [N]} \Xi_i) \ | \ \mathbb{E}_{\mathbb{P}}[f_{i,l}(\tilde{\xi}_i)] \in [\underline{\gamma}_{i,l},\overline{\gamma}_{i,l}], \forall l \in [L_i], \forall i \in [N],\\
\quad \quad \quad \quad \quad  \mathbb{E}[f_{i,j,l}(\tilde{\xi}_i,\tilde{\xi}_j)] \leq \gamma_{i,j,l}, \forall l \in [L_{i,j}], \forall (i,j) \in [N]_2\big\},
\end{array}
\end{equation}
where the bivariate functions are submodular. As discussed in Section \ref{sec3}, both marginal information and bivariate dependence information is modeled using ${\cal P}_{\text{uni}} \cap {\cal P}_{\text{bi-sub}}$.

In the following theorem, we show that problem \eqref{droiner00} is solvable efficiently as a compact linear program under appropriate assumptions on the objective function. 
\begin{theorem} \label{thm:affine}
Suppose the off-diagonal entries of the matrices $\mb{A}_k$ are nonnegative for each $k \in [K]$. When $\Xi_i$ is discrete and finite for each $i \in [N]$, $\rho^*$ in \eqref{droiner00} is given by the optimal value of the polynomial sized linear program:
\begin{equation}\label{primalmdm1}
\begin{array}{rllll}
\displaystyle \max_{\mbs{\lambda}} & \displaystyle \sum_{k \in [K]}\left(\sum_{(i,j) \in [N]_2}\sum_{\xi_i \in {\Xi}_i}\sum_{\xi_j \in {\Xi}_j}  2A_{i,j,k}\xi_i\xi_j\lambda_{i,j,k}(\xi_i,\xi_j)+\sum_{i \in [N]}\sum_{\xi_i \in {\Xi}_i} \left(A_{i,i,k}\xi_i^2+b_{i,k}\xi_i\right)\lambda_{i,k}(\xi_i)+c_k\lambda_k \right)& \\
\textrm{s.t.} & \displaystyle \sum_{k \in [K]} \lambda_k = 1,\\
& \displaystyle  \sum_{\xi_i \in \Xi_i}\lambda_{i,k}(\xi_i) = \lambda_k, \quad \quad \quad \quad \quad  \quad \quad \quad  \quad \quad \quad   \quad \quad \quad 
 \forall i \in [N], \forall k \in [K],\\
& \displaystyle \underline{\gamma}_{i,l} \leq \sum_{k \in [K]}\sum_{\xi_i \in \Xi_i}f_{i,l}(\xi_i)\lambda_{i,k}(\xi_i) \leq \overline{\gamma}_{i,l},  \quad \quad \quad \quad \quad \quad \forall l \in [L_i], \forall i \in [N],\\
& \displaystyle \sum_{\xi_i \in \Xi_i}\lambda_{i,j,k}(\xi_i,\xi_j) = \lambda_{j,k}(\xi_j), \quad \quad \quad  \quad \quad \quad \quad  \quad \quad \quad   \forall \xi_j \in \Xi_j, \forall (i,j)  \in [N]_2, \forall k \in [K],\\
& \displaystyle \sum_{\xi_j \in \Xi_j}\lambda_{i,j,k}(\xi_i,\xi_j) = \lambda_{i,k}(\xi_i),  \quad \quad \quad  \quad \quad \quad \quad  \quad \quad \quad\forall \xi_i \in \Xi_i, \forall (i,j)  \in [N]_2, \forall k \in [K],\\
& \displaystyle  \sum_{k \in [K]}\sum_{\xi_i \in \Xi_i}\sum_{\xi_j \in \Xi_j}f_{i,j,l}({\xi}_{i},{\xi}_{j})\lambda_{i,j,k}(\xi_i,\xi_j) \leq \gamma_{i,j,l},  \ \   \quad 
 \quad  \forall l \in [L_{i,j}], \forall (i,j) \in [N]_2,\\
& \displaystyle  \lambda_k \geq 0,  \ \ \quad \quad \quad \quad \quad \quad\quad \quad \quad  \quad \quad \quad \quad \quad \quad  \quad \quad \quad 
\forall k \in [K],\\
& \displaystyle  \lambda_{i,k}(\xi_i) \geq 0,\ \ \quad \quad \quad \quad\quad \quad \quad  \quad \quad \quad\quad \quad \quad  \quad \quad \quad \forall \xi_i \in \Xi_i, \forall i \in [N], \forall k \in [K],\\
& \displaystyle  \lambda_{i,j,k}(\xi_i,\xi_j) \geq 0,  \quad \quad \quad \quad\quad \quad \quad  \quad \quad \quad \quad \quad \quad \quad \quad \forall \xi_i \in \Xi_i, \forall \xi_j \in \Xi_j, \forall (i,j) \in [N]_2, \forall k \in [K].
\end{array}
\end{equation}
\end{theorem}
\begin{proof}
See Appendix.
\end{proof}

We make a few remarks about Theorem \ref{thm:affine} and its implications next.
\texitem{(a)} The proof is constructive where the extremal distribution is given by a mixture of comonotonic random vectors (see Figures \ref{fig:tail2a} and \ref{fig:tail2b} for an illustration of the construction). However the extremal distribution itself need not be comonotonic; also see Example 7.
\begin{figure}[!htbp]
\centering
\subfigure[$k = 1$] {
\begin{tikzpicture}[scale=0.6]
\draw(0.64*4,1)node[right]{{$F_{2|k}^{-1}(U)$}};
\draw(0.64*4,4)node[left]{{$F_{1|k}^{-1}(U)$}};
 \draw[-] (0,0) -- (4,0)node[below]{{1}} -- (4,0)node[right]{$U$} ;
  \draw[-] (0,0)node[below]{{$0$}} -- (0,4) node[left]{{$F^{-1}(U)$}};
 \draw[dashed] (0,0) -- (0.36*4,0) -- (0.36*4,1) --  (0.64*4,1)  -- (0.64*4,2) -- (0.84*4,2) -- (0.84*4,3)-- (0.96*4,3) -- (0.96*4,4) -- (1*4,4);
 \draw[solid] (0,0) -- (0.04*4,0) -- (0.04*4,1) --  (0.16*4,1)  -- (0.16*4,2) -- (0.36*4,2) -- (0.36*4,3)-- (0.64*4,3) -- (0.64*4,4) -- (1*4,4);
 \end{tikzpicture}}
\subfigure[$k = 2$] {\begin{tikzpicture}[scale=0.6]
\draw(0.64*4,1)node[right]{{$F_{1|k}^{-1}(U)$}};
\draw(0.64*4,4)node[left]{{$F_{2|k}^{-1}(U)$}};
  \draw[-] (0,0) -- (4,0)node[below]{{$1$}} -- (4,0)node[right]{$U$} ;
   \draw[-] (0,0)node[below]{{$0$}} -- (0,4) node[left]{{$F^{-1}(U)$}};
 \draw[dashed] (0,0) -- (0.04*4,0) -- (0.04*4,1) --  (0.16*4,1)  -- (0.16*4,2) -- (0.36*4,2) -- (0.36*4,3)-- (0.64*4,3) -- (0.64*4,4) -- (1*4,4);
  \draw[solid] (0,0) -- (0.36*4,0) -- (0.36*4,1) --  (0.64*4,1)  -- (0.64*4,2) -- (0.84*4,2) -- (0.84*4,3)-- (0.96*4,3) -- (0.96*4,4) -- (1*4,4);
\end{tikzpicture}
}
\caption{Comonotonic construction for the conditional distributions in step (ii) for $N = 2$ and $K = 2$. Here the solid line indicates $\tilde{\xi}_1 \sim F_{1|k}^{-1}(U)$ and the dashed line indicates $\tilde{\xi}_2 \sim F_{2|k}^{-1}(U)$ for $k = 1$ (left figure) and $k = 2$ (right figure) where $F_{i|k}$ is the conditional marginal distribution of $\tilde{\xi}_i$ for index $k$ being optimal.} \label{fig:tail2a}
\end{figure}
\begin{figure}[!htbp]
\centering
\subfigure[$k = 1$] {
\begin{tikzpicture}[scale=0.6]
\draw[step=1cm,dotted] (0,0) grid (4,4);
\node[below] at (4,0) {$\xi_1$};
\node[left] at (0,4) {$\xi_2$};
\filldraw[black](0,0)circle[radius=2pt]; 
\filldraw[black](1,0)circle[radius=2pt]; 
\filldraw[black](2,0)circle[radius=2pt]; 
\filldraw[black](3,1)circle[radius=2pt]; 
\filldraw[black](4,2)circle[radius=2pt];  
\filldraw[black](4,3)circle[radius=2pt]; 
\filldraw[black](4,4)circle[radius=2pt]; 
\node at (5,2) {$+$};
\end{tikzpicture}
}
\subfigure[$k = 2$] {
\begin{tikzpicture}[scale=0.6]
\draw[step=1cm,dotted] (0,0) grid (4,4);
\node[below] at (4,0) {$\xi_1$};
\node[left] at (0,4) {$\xi_2$};
\filldraw[black](0,0)circle[radius=2pt]; 
\filldraw[black](0,1)circle[radius=2pt]; 
\filldraw[black](0,2)circle[radius=2pt]; 
\filldraw[black](1,3)circle[radius=2pt]; 
\filldraw[black](2,4)circle[radius=2pt];  
\filldraw[black](3,4)circle[radius=2pt]; 
\filldraw[black](4,4)circle[radius=2pt]; 
\node at (5,2) {$\Longrightarrow$};
\end{tikzpicture}
}
\subfigure[Overall] {
\begin{tikzpicture}[scale=0.6]
\draw[step=1cm,dotted] (0,0) grid (4,4);
\node[below] at (4,0) {$\xi_1$};
\node[left] at (0,4) {$\xi_2$};
\filldraw[black](0,0)circle[radius=2.5pt]; 
\filldraw[black](1,0)circle[radius=2pt]; 
\filldraw[black](2,0)circle[radius=2pt]; 
\filldraw[black](3,1)circle[radius=2pt]; 
\filldraw[black](4,2)circle[radius=2pt];  
\filldraw[black](4,3)circle[radius=2pt]; 
\filldraw[black](4,4)circle[radius=2.5pt]; 
\filldraw[black](0,1)circle[radius=2pt]; 
\filldraw[black](0,2)circle[radius=2pt]; 
\filldraw[black](1,3)circle[radius=2pt]; 
\filldraw[black](2,4)circle[radius=2pt];  
\filldraw[black](3,4)circle[radius=2pt]; 
\end{tikzpicture}
}
\caption{Subfigures (a) and (b) display the support of the conditional bivariate distributions for $k = 1$ and $k = 2$ while (c) shows that overall support of the extremal bivariate distribution using the weighted probabilities $\lambda_1^*$ and $\lambda_2^*$. While the conditional bivariate distributions in (a) and (b) are comonotonic, the final bivariate distribution in (c) is not comonotonic but positively dependent here.} \label{fig:tail2b}
\end{figure}
\texitem{(b)} The proof shows that the ambiguity set  ${\cal P}_{\text{uni}} \cap {\cal P}_\text{bi-sub}$ is nonempty if and only if there exists scalars $p_{i}(\xi_i)$ for all $\xi_i \in \Xi_i$, $i \in [N]$ and $p_{i,j}(\xi_i,\xi_j)$ for all $\xi_i \in \Xi_i,\xi_j \in \Xi_j$, $(i,j) \in [N]_2$ such that the following linear conditions are jointly satisfied:
\begin{equation*}
\begin{array}{rllll}
& \displaystyle  \sum_{\xi_i \in \Xi_i}p_{i}(\xi_i) = 1, & 
 \forall i \in [N], \\
& \displaystyle \underline{\gamma}_{i,l} \leq \sum_{\xi_i \in \Xi_i}f_{i,l}(\xi_i)p_{i}(\xi_i) \leq \overline{\gamma}_{i,l},& \forall l \in [L_i], \\
& \displaystyle \sum_{\xi_i \in \Xi_i}p_{i,j}(\xi_i,\xi_j) = p_j(\xi_j), &  \forall \xi_j \in \Xi_j, \forall (i,j)  \in [N]_2,\\
& \displaystyle \sum_{\xi_j \in \Xi_j}p_{i,j}(\xi_i,\xi_j) = p_{i}(\xi_i),  & \forall \xi_i \in \Xi_i, \forall (i,j)  \in [N]_2,\\
& \displaystyle  \sum_{\xi_i \in \Xi_i}\sum_{\xi_j \in \Xi_j}f_{i,j,l}({\xi}_{i},{\xi}_{j})p_{i,j}(\xi_i,\xi_j) \leq \gamma_{i,j,l},  &  \forall l \in [L_{i,j}], \forall (i,j) \in [N]_2,\\
& \displaystyle  p_{i}(\xi_i) \geq 0, &\forall \xi_i \in \Xi_i, \forall i \in [N],\\
& \displaystyle  p_{i,j}(\xi_i,\xi_j) \geq 0,  &\forall \xi_i \in \Xi_i, \forall \xi_j \in \Xi_j, \forall (i,j) \in [N]_2.
\end{array}
\end{equation*}
This is verifiable in polynomial time.
\texitem{(c)} The size of the linear program in Theorem \ref{thm:affine} grows quadratically in $N$ (number of random variables) and $\max_{i \in [N]} |\Xi_i|$ (maximum number of values that a random variable can take) and linearly in $K$ (number of pieces defining the objective function)   and $\max_{i \in [N]} L_i$ and $\max_{(i,j) \in [N]_2} L_{i,j}$ (number of constraints in the submodular ambiguity set). The proof builds on the proof of Theorem 2.3.1 in \cite{karthikbook}. However the linear program therein is developed only for ${\cal P}_{\text{uni}}$. The key novelty in Theorem \ref{thm:affine} is that by introducing additional decision variables for modeling bivariate probabilities and incorporating appropriate constraints, we show that it is possible to compute the sharp bound for ${\cal P}_{\text{uni}} \cap {\cal P}_{\text{bi-sub}}$ in polynomial time with linear programming. 
\texitem{(d)} It is straightforward to build on Theorem \ref{thm:affine} to solve distributionally robust optimization problems. Specifically, consider:
\begin{equation} \label{dro11}
\begin{array}{rlll}
\displaystyle \inf_{\mbs{x} \in {\cal X}} \sup_{\mathbb{P} \in {\cal P}_{\text{uni}} \cap {\cal P}_\text{bi-sub}}\mathbb{E}_{\mathbb{P}}\left[\max_{k \in [K]}\left(\tilde{\mb{\xi}}'\mb{A}_k(\mb{x})\tilde{\mb{\xi}} + \mb{b}_k'(\mb{x})\tilde{\mb{\xi}} + c_{k}(\mb{x})\right)\right],
\end{array}
\end{equation}
where for each $k \in [K]$, the matrix $\mb{A}_k(\mb{x})$, the vector $\mb{b}_k(\mb{x})$ and the scalar $c_k(\mb{x})$ have an affine dependence on $\mb{x}$. Assume the set ${\cal X}$ is given by the polyhedron $\{\mb{x} \ | \ \mb{D}\mb{x} \geq \mb{d}\}$ and the off-diagonal entries of the matrix $\mb{A}_k(\mb{x})$ are nonnegative for each $k \in [K]$ and $\mb{x} \in {\cal X}$. Then by using the dual of the linear program in \eqref{primalmdm1}, we can reformulate \eqref{dro11} as a linear program. Such a technique is commonly used in robust and distributionally robust optimization.(see \cite{bentalbook,hertog}).

\section{Continuous uncertainty}
\label{sec5}
In this section, we discuss models for the distributionally robust optimization problem with a submodular ambiguity set when the support of the uncertainty is allowed to be continuous. 
\subsection{Pseudo-polynomial time solvability with discretization} \label{sec5a}
We first focus on solving the problem \eqref{droiner} and make the following assumptions on the ambiguity set and the cost function:\\
\texttt{Assumption (A1'):} The set $\Xi_i =[0,B_i]\subset \mathbb{R}$ is a finite interval for each $i \in [N]$.\\
\texttt{Assumption (A2'):} For each $l \in [L]$, $f_l: \prod_{i \in [N]} \Xi_i \rightarrow \mathbb{R}$ is a $M$-Lipschitz submodular function with respect to the infinity norm with a polynomial time evaluation oracle. \\
\texttt{Assumption (A3'):} The ambiguity set ${\cal P}_{\text{sub}}$ is nonempty with $\sum_{r \in [R]} p_r \delta_{\mbs{\xi}_r} \in {\cal P}_{\text{sub}}$.\\
\texttt{Assumption (A4'):} The cost function is given by $g({\mb{\xi}}) := \max_{k \in [K]}g_k(\mb{\xi})$ where for each $k \in [K]$, the function $g_k: \prod_{i \in [N]} \Xi_i \rightarrow \mathbb{R}$ is a $M$-Lipschitz supermodular function with respect to the infinity norm with a polynomial time evaluation oracle.\\
This brings us to the following theorem.
\begin{theorem} \label{thm:submodcon}
Suppose assumptions (A1')-(A4') hold. Then for every $\epsilon > 0$, one can find a value $\rho_{\epsilon}^* \geq \rho^*-\epsilon$ in time polynomial in $N$, $L$, $K$, $\max_{i \in [N]}B_i$, $R$, $M$, $\frac{1}{\epsilon}$ and the evaluation time of the oracles such that there exists a distribution $\mathbb{P}_{\epsilon}^* \in {\cal P}(\prod_{i \in [N]}\Xi_i)$ with $\rho_{\epsilon}^* = \mathbb{E}_{\mathbb{P}_{\epsilon}^*}[g(\tilde{\mb{\xi}})]$ and $\mathbb{E}_{\mathbb{P}_{\epsilon}^*}[f_l(\tilde{\mb{\xi}})] \leq \gamma_l+\epsilon$ for all $l \in [L]$.
\end{theorem}
\begin{proof}
See Appendix.
\end{proof}

\texitem{(a)}
 The proof technique in Theorem \ref{thm:submodcon} builds on the discretization approach first proposed in \cite{bach} for deterministic Lipschitz continuous submodular function minimization. The usefulness of the theorem lies in that even when the uncertainty is defined on a continuous support, it is possible to obtain pseudo-polynomial time $\epsilon$-additive guarantees. However much stronger polynomial time solution guarantees can be obtained for continuous uncertainty with a convex ambiguity set. On the other hand, the theorem shows that it is possible to incorporate information on the uncertainty that might not be possible with a convex ambiguity set while still being solvable in pseudo-polynomial time. As an illustration, consider the
ambiguity set in Example 5. This is not a convex ambiguity set but is a submodular ambiguity set and Theorem \ref{thm:submodcon} applies. To the best of our knowledge, such a result is not covered by existing work in
distributionally robust optimization
\texitem{(b)} When the submodular ambiguity set is defined using only univariate and bivariate information and the support of each random variable is an interval, we can use discretization for the marginal support and solve the linear program provided in Theorem \ref{thm:affine} to obtain the additive approximation. The size of the linear program in this case is pseudo-polynomial in the size of the input. It is also straightforward to extend these results to the distributionally robust optimization problem through duality. In the next section, we provide a structured submodular ambiguity set defined using moment information with continuous uncertainty, where one can obtain stronger polynomial time solution guarantees using semidefinite programming.

\subsection{Polynomial sized semidefinite program with continuous uncertainty}  \label{sec3ba}
In this section, we propose a moment based submodular ambiguity set for nonnegative random vectors. Given a vector $\mb{\mu} \in \mathbb{R}^N_+$ and a symmetric matrix $\mb{\Sigma} \in \mathbb{S}^N$ (not necessarily positive semidefinite) with nonnegative entries along the diagonal, define a cross moment ambiguity set ${\cal Q}_{\text{cm}}$ as follows:
\begin{equation} \label{eq:bimom}
\begin{array}{lll}
{\cal Q}_{\text{cm}} :=  \{\mathbb{P} \in {\cal P}(\mathbb{R}^N_+) \ | \ \mathbb{E}_{\mathbb{P}}[\tilde{\mb{\xi}}] = \mb{\mu}, \mathbb{E}_{\mathbb{P}}[\mbox{diag}(\tilde{\mb{\xi}}\tilde{\mb{\xi}}')] = \mbox{diag}(\mb{\mu}\mb{\mu}'+\mb{\Sigma}), \mathbb{E}_{\mathbb{P}}[(\tilde{\mb{\xi}}-\mb{\mu})(\tilde{\mb{\xi}}-\mb{\mu})')] \geq \mb{\Sigma} \big\}.
 \end{array}
 \end{equation}
 The mean of the random variable $\tilde{\xi}_i$ is fixed to $\mu_i$ and the variance is fixed to $\Sigma_{i,i}$ for each $i \in [N]$. On the other hand, the covariance for each pair of random variables is bounded from below with $\mbox{Cov}(\tilde{\xi}_i,\tilde{\xi}_j) \geq \Sigma_{i,j}$. The set \eqref{eq:bimom} is obtained by intersecting an univariate and a bivariate submodular ambiguity set. Now, consider the problem of computing the tightest upper bound on the expected value of a piecewise quadratic function of the random vector for the cross moment ambiguity set given by:
\begin{equation} \label{droinerabc}
\begin{array}{rlll}
\displaystyle \rho^* := \sup_{\mathbb{P} \in {\cal Q}_{\text{cm}}}\mathbb{E}_{\mathbb{P}}\left[\max_{k \in [K]}\left(\tilde{\mb{\xi}}'\mb{A}_k\tilde{\mb{\xi}}+\mb{b}_k'\tilde{\mb{\xi}} + b_{k}\right)\right].
\end{array}
\end{equation}
 This brings us to the following theorem on computing the moment bound in \eqref{droinerabc} in polynomial time using semidefinite optimization.
 \begin{theorem} \label{thm:affinemomsdp}
Suppose the off diagonal entries of the matrices $\mb{A}_k$ are nonnegative for each $k \in [K]$.
Then $\rho^*$ in \eqref{droinerabc} is given by the optimal value of the polynomial sized semidefinite program:
\begin{equation}\label{primalmdm1momsdpa}
\begin{array}{rllll}
\displaystyle \max_{\mbs{\lambda}} & \displaystyle \sum_{k \in [K]}\left(\sum_{(i,j) \in [N]_2} 2A_{i,j,k}\lambda_{i,j,k}+\sum_{i \in [N]}\left(A_{i,k}\lambda_{i,i,k}+ b_{i,k}\right)\lambda_{i,k}+c_k\lambda_k\right)& \\
\mbox{s.t.} & \displaystyle \sum_{k \in [K]} \lambda_k = 1,\\
& \displaystyle  \sum_{k \in [K]}{\lambda}_{i,k} = \mu_i, \quad \quad \quad \quad \quad \quad \quad \forall i \in [N],\\
& \displaystyle  \sum_{k \in [K]}{\lambda}_{i,i,k} = \mu_i^2+\Sigma_{i,i},  \quad \quad \quad \quad\forall i \in [N], \\
& \displaystyle \sum_{k \in [K]}{\lambda}_{i,j,k} \geq \mu_i\mu_j+\Sigma_{i,j}, \quad \quad \quad \forall (i,j) \in [N]_2, \\
&   \displaystyle \begin{pmatrix}
 \lambda_k & \lambda_{i,k}  & \lambda_{j,k}\\
  \lambda_{i,k} & \lambda_{i,i,k}  & \lambda_{i,j,k}\\
   \lambda_{j,k} & \lambda_{i,i,k}  & \lambda_{j,j,k}
    \end{pmatrix} \succeq 0,\  \quad \quad \forall (i,j) \in [N]_2,\forall k \in [K],\\
    &\displaystyle \lambda_{i,k} \geq 0, \ \ \  \quad \quad \quad \quad \quad \quad \quad \quad \quad \forall i  \in [N],\forall k \in [K].
\end{array}
\end{equation}
\end{theorem}
\begin{proof}
See Appendix.
\end{proof}

We make a few remarks about Theorem \ref{thm:affinemomsdp} and its implications next.
\texitem{(a)} In comparison to Theorem \ref{thm:affine} where probabilities are used as decision variables, the formulation above makes use of moments as decision variables. However here the support of the random vector is allowed to lie in $\mathbb{R}^N_+$ and hence the corresponding moment cone is not closed. Specifically, define the moment cone ${\cal M}$ as follows:
\begin{equation*} \label{eq:mom}
{\cal M} :=
\left\{
    z_0(1,\mb{z},\mb{Z}) \ | \ \exists \mathbb{P} \in {\cal P}(\mathbb{R}^N_+) \mbox{ s.t. }\mathbb{E}_{\mathbb{P}}[\tilde{\mb{\xi}}] = \mb{z}, \mathbb{E}_{\mathbb{P}}[\mbox{diag}(\tilde{\mb{\xi}}\tilde{\mb{\xi}}')] = \mbox{diag}(\mb{Z}), \mathbb{E}_{\mathbb{P}}[\tilde{\mb{\xi}}\tilde{\mb{\xi}}'] \geq \mb{Z},  z_0 \geq 0
\right\}.
\end{equation*}
The proof shows that $(1,\mb{\mu},\mb{\mu}\mb{\mu}'+\mb{\Sigma}) \in \mbox{cl}$(${\cal M}$) if and only if there exists scalars $m_{i,j}$ for $(i,j) \in [N]_2$ such that the following conditions are jointly satisfied:
\begin{equation*}
\begin{array}{rllll}
& \displaystyle {m}_{i,j} \geq \mu_i\mu_j+\Sigma_{i,j}, &\forall (i,j) \in [N]_2, \\
&   \displaystyle \begin{pmatrix}
 1 & \mu_{i}  & \mu_{j}\\
  \mu_{i} & \mu_i^2+\Sigma_{i,i}  & m_{i,j}\\
   \mu_{j} & m_{i,j}  & \mu_j^2+\Sigma_{j,j}
    \end{pmatrix} \succeq 0, &\forall (i,j) \in [N]_2,\\
    &\displaystyle \mu_{i} \geq 0, & \forall i  \in [N].
\end{array}
\end{equation*}
As an example, let $N = 2$, $\mu_1  = \mu_2= 0$, $\Sigma_{1,1} = \Sigma_{2,2}= \Sigma_{1,2}= 1$. This is a point which lies only in the closure of the moment cone. To see this, observe that the moments are attained by the limit of a sequence of probability measures of the form below where $\epsilon \downarrow 0$: 
\begin{equation*}
\begin{array}{llllll}
\displaystyle \mathbb{P}(\tilde{\xi}_1 = 0,\tilde{\xi}_2 =0)   & = & 1-\epsilon, \\
\mathbb{P}\left(\tilde{\xi}_1 = \frac{1}{\sqrt{\epsilon}},\tilde{\xi}_2 = \frac{1}{\sqrt{\epsilon}}\right)   &= & \epsilon.
\end{array}
\end{equation*}
 In the proof, we provide two constructions - the first construction works in the simpler case when the bound is attained by an extremal distribution and in the second and more general case, we construct a sequence of distributions that in the limit attains the bound and satisfies the moment conditions. The result in Theorem \ref{thm:affinemomsdp} extends in a straightforward manner to support in $\mathbb{R}^N$ where the last set of nonnegativity conditions in \eqref{primalmdm1momsdpa} is dropped.
\texitem{(b)} The ambiguity set ${\cal Q}_{\text{cm}}$ in \eqref{eq:bimom} can be equivalently rewritten as:
\begin{equation*} \label{eq:bimom0}
\begin{array}{lll}
{\cal Q}_{\text{cm}} =  \{\mathbb{P} \in {\cal P}(\mathbb{R}^N_+) \ | \ \mathbb{E}_{\mathbb{P}}[\tilde{{\xi}}_i] = {\mu}_2, \forall i \in [N], \mathbb{E}_{\mathbb{P}}[\tilde{{\xi}}_i^2] = \mu_i^2+\Sigma_{i,i}, \forall i \in [N],\\
\quad \quad \quad \quad \quad \quad \quad \quad \quad \mathbb{E}_{\mathbb{P}}[(\tilde{{\xi}}_i-\tilde{{\xi}}_j)^2] \leq \mu_i+\Sigma_{i,i}+\mu_j+\Sigma_{j,j}-2\mu_i\mu_j-2\Sigma_{i,j}, \forall (i,j) \in [N]_2 \big\},
 \end{array}
 \end{equation*}
 where upper bounds on pairwise Wasserstein distances are given. While this is a submodular ambiguity set, it is not an example of a convex ambiguity set since the second moments of $\tilde{\xi}_i$ are precisely specified. 
\texitem{(c)} The semidefinite program in \eqref{primalmdm1momsdpa} is closely related to the semidefinite program that has been developed for the cross moment ambiguity set ${\cal P}_{\text{cm}}$ discussed in Section \ref{sec2}. When the support is $\Xi = \mathbb{R}^N_+$ with $\mb{\mu} \in \mathbb{R}^N_+$ and $\mb{\Sigma} \succeq 0$, it was shown in \cite{Delage} that:
\begin{equation}\label{primalmdm1momsdpabb}
\begin{array}{rllll}
\displaystyle \sup_{\mathbb{P} \in {\cal P}_{\text{cm}}}\mathbb{E}_{\mathbb{P}}\left[\max_{k \in [K]}\left(\tilde{\mb{\xi}}'\mb{A}_k\tilde{\mb{\xi}}+\mb{b}_k'\tilde{\mb{\xi}} + c_{k}\right)\right] =
\displaystyle \max_{\mbs{\lambda},\mbs{\Lambda}} & \displaystyle \sum_{k \in [K]} \mb{A}_k \bullet \mb{\Lambda}_k+\sum_{k \in [K]}\mb{b}_k'\mb{\lambda}_k +\sum_{k \in [K]}c_k\lambda_k& \\
\mbox{s.t.} & \displaystyle \sum_{k \in [K]} \lambda_k = 1,\\
& \displaystyle \sum_{k \in [K]}\mb{\lambda}_{k} = \mb{\mu},\\
& \displaystyle \sum_{k \in [K]}\mb{\Lambda}_{k} \preceq \mb{\mu}\mb{\mu}'+\mb{\Sigma},\\
& \displaystyle \mb{\lambda}_k \geq 0,\quad \quad  \quad \quad  \quad \quad \forall k  \in [K],\\
&   \displaystyle \begin{pmatrix}
 \lambda_k & \mb{\lambda}_k'\\
  \mb{\lambda}_k &  \mb{\Lambda}_k
    \end{pmatrix} \succeq 0, \quad \quad \forall k \in [K],
\end{array}
\end{equation}
where the matrices $\mb{A}_k$ are assumed to be negative semidefinite. In comparing \eqref{primalmdm1momsdpa} with \eqref{primalmdm1momsdpabb}, we observe that the former formulation uses $O(N^2K)$ positive semidefinite matrices of size $3 \times 3$ while the latter formulation uses $O(K)$ positive semidefinite matrices of size $(N+1) \times (N+1)$. The formulation \eqref{primalmdm1momsdpa} incorporates lower bounds on the covariances term by term and for tightness of the formulation, requires the matrices $\mb{A}_k$ to satisfy sign restrictions on the off-diagonal entries. On the other hand, formulation \eqref{primalmdm1momsdpabb} specifies upper bounds in a positive semidefinite order and requires the matrices $\mb{A}_k$ to be negative semidefinite for tightness. The related problem of testing if there exists a nonnegative random vector with a given mean vector $\mb{\mu}$ and covariance matrix $\mb{\Sigma}$ and optimizing over this set of distributions is known to be NP-hard (see \cite{dickinson}). 

\section{Numerical Examples} \label{sec6}

\subsection{Project network} We discuss a flexible approach to model activity duration distributions in PERT networks while incorporating lower bounds on the correlation among the activities.
Consider the network in Figure \ref{fig:pert}, first proposed in \cite{slyke}. This network consists of seven activities denoted by the arcs and five paths where four parallel arcs are present in the top part of the network. The start of the project is denoted by node 1 and the end by node 4. The length of an arc denotes the duration of the corresponding activity and the length of the longest path from node 1 to node 4 provides the project completion time. We are interested in the case where the activity durations are random.
\begin{figure}[h]
\centering
\begin{tikzpicture}
\tikzset{vertex/.style = {shape=circle,draw,minimum size=0.2em}}
\tikzset{edge/.style = {->,> = latex}}
\node[vertex] (C1) at  (0,0) {$1$};
\node[vertex] (C2) at  (2.5,2) {$3$};
\node[vertex] (C3) at  (2.5,-2) {$2$};
\node[vertex] (C4) at  (5,0) {$4$};
\draw[edge]  (C1) to (C2) ;
\draw[edge] (C1) to (C3) ;
\draw[edge] (C3) to (C4) ;
\draw[edge] (C2) to [out=30,in=70] (C4) ;
\draw[edge] (C2) to [out=-10,in=100] (C4) ;
\draw[edge] (C2) to [out=-60,in=170] (C4) ;
\draw[edge] (C2) to [out=-90,in=210] (C4) ;
\end{tikzpicture}
\caption{PERT network with seven activities and five paths from start node $1$ to end $4$.} \label{fig:pert}
\end{figure}
We compute the worst-case expected project duration where information on the marginal distribution of the activity durations is available along with correlation information. Four types of marginal distributions are considered in the numerical study:
\begin{enumerate}
\item Three point distribution: Duration of activity $(1,2)$ is supported on $\{5,10,15\}$ with probabilities $(1/6,4/6,1/6)$ while all other activities have durations with support in $\{4,9,14\}$ with probabilities $(1/6,4/6,1/6)$.
\item Uniform distribution: Duration of activity $(1,2)$ is uniform in the range $[5,15]$ while duration of all other activities is uniform in the range  $[4,14]$.  
\item Beta distribution: We fit continuous PERT beta distributions to match the range and mode of the three point distributions with variance fixed to 1/6 of the range. 
\item Triangular distribution: We fit continuous triangular distributions to match the range and mode of the three point distributions. 
\end{enumerate}
\noindent To deal with the uniform, beta and triangular distributions, we discretize the marginal distributions uniformly in the range. The number of marginal support points for each activity duration is identical and varied in the set $\{10,20,30,40,50\}$. We compute the worst-case expected project duration using the linear program provided in \eqref{primalmdm1} where the lower bound on the correlations between all pairs of activities is varied in $\{-1,0,0.2,0.4,0.6,0.8,1\}$. When the lower bound is $-1$, the results correspond to worst-case estimates for the univariate ambiguity set where only the marginals are given. When the lower bound on the correlation is $1$, the results correspond to the comonotonic distribution which in this network gives an expected duration of 19 and a critical path of activities $(1,2)$ and $(2,4)$. Solving the linear programs for the three point distribution gives the worst-case expected project duration values of $\{24.17, 24.17, 23.82, 23.18, 21.80, 20.39, 19.0\}$ and the criticality of activity $(1,2)$ is given by $\{0.67, 0.67, 0.68, 0.79, 0.53,0.27, 0\}$ for the different lower bounds on correlations. In this network, with stochastic activity durations, activity $(1,3)$ is more critical than activity $(1,2)$ for lower bounds on correlations of $\{-1,0,0.2,0.4, 0.6\}$ due to the presence of the parallel activities in the top part of the network while for correlations of $\{0.8,1\}$, activity $(1,2)$ is more critical than $(1,3)$ due to the positive shift in the marginal distribution of activity $(1,2)$ as compared to $(1,3)$. The results for the other three marginal distributions are provided in Tables \ref{table:uniform}, \ref{table:beta} and \ref{table:triangular} where the worst-case expected project duration and the criticality of activity $(1,3)$
 is reported. The tables illustrate that correlation information on activity durations, when available, is important to incorporate along with the marginal information as the estimates are sensitive to this information. As the lower bounds on the correlation increases (namely activities are allowed to be more positively correlated), the worst-case project duration decreases in this example and the criticality of activity $(1,3)$ also decreases. The worst-case estimates for the uniform distribution are higher as compared to the beta and triangular distributions, indicating the flexibility in capturing marginal information in obtaining project duration estimates. Lastly, the results indicate that as the discretization increases, the marginal benefits obtained from solving the larger linear programs decreases. 

\begin{table}[!htbp]
\centering
\begin{tabular}{|l|l|l|l|l|l|l|}
\hline
Correlation & Support points& 10 & 20 & 30 & 40 & 50 \\\hline
-1 & Duration & 24.50 & 24.55 & 24.55 &  24.55 &  24.55 \\
& Criticality &0.60  & 0.60  & 0.59  & 0.59  &0.58 \\ \hline
0 & Duration &24.03 &  24.10 &  24.10 &  24.10 & 24.11 \\
& Criticality &0.68 &  0.65 &  0.67 & 0.65 & 0.64\\\hline
0.2 & Duration &23.80 &  23.86 &  23.86 &  23.86 & 23.86 \\
& Criticality & 0.70 & 0.65  & 0.66 & 0.66 &  0.66\\\hline
0.4 & Duration &23.44 &  23.50 &  23.51 &  23.52 &23.52 \\
& Criticality &0.60 &  0.62 &  0.63 & 0.63 & 0.63\\\hline
0. & Duration &22.57 & 22.63 &  22.64 & 22.64 & 22.64 \\
& Criticality &0.58 &  0.58 &  0.57 & 0.58 &  0.58 \\\hline
0.8 & Duration & 21.36 & 21.42 &  21.43 &  21.44 & 21.44 \\
& Criticality &0.55 & 0.55  &0.54 &  0.54 & 0.54 \\\hline
1 & Duration &19.0 &  19.0 &  19.0  &19.0 & 19.0\\
& Criticality &0 & 0&  0& 0& 0 \\
\hline
\end{tabular}
\caption{Uniform activity durations: Worst-case expected project duration and criticality of activity (1,3) with given lower bound on correlation and number of support points.}  \label{table:uniform}
\end{table}

\begin{table}[!htbp]
\centering
\begin{tabular}{|l|l|l|l|l|l|l|}
\hline
Correlation & Support points& 10 & 20 & 30 & 40 & 50 \\\hline
-1 & Duration & 21.17 &  21.15 &  21.13 &  21.13 & 21.14 \\
& Criticality &0.50 &  0.51 &  0.58 &  0.57 & 0.56\\ \hline
0 & Duration &20.97 & 21.00 & 21.00 &   21.00 & 21.00  \\
& Criticality & 0.56 & 0.55 & 0.58 &   0.56 &0.58\\\hline
0.2 & Duration & 20.85&   20.89 &   20.89 &  20.889 & 20.89 \\
& Criticality & 0.56 & 0.54 &0.58 &0.56 &0.58\\\hline
0.4 & Duration &20.74 &   20.75 &  20.75 &   20.75 &20.75\\
& Criticality &0.56 &   0.54 &  0.56 &   0.55 &   0.56\\\hline
0.6 & Duration &20.37 & 20.38 &  20.38 &  20.38 &  20.39 \\
& Criticality & 0.53 &   0.49 &  0.48 &   0.49 &  0.48\\\hline
0.8 & Duration & 19.756 &   19.85 & 19.856 &   19.86 &  19.86 \\
& Criticality &0.35 &  0.38 &  0.39 &  0.38&  0.38 \\\hline
1 & Duration &19.0 &  19.0 &  19.0  &19.0 & 19.0\\
& Criticality &0 & 0&  0& 0& 0 \\
\hline
\end{tabular}
\caption{Beta activity durations: Worst-case expected project duration and criticality of activity (1,3) with given lower bound on correlation and number of support points.}  \label{table:beta}
\end{table}

\begin{table}[!htbp]
\centering
\begin{tabular}{|l|l|l|l|l|l|l|}
\hline
Correlation & Support points& 10 & 20 & 30 & 40 & 50 \\\hline
-1 & Duration & 22.70 &  21.71 &   21.72 &   21.72 &  21.72 \\
& Criticality &0.61 &   0.56 & 0.56 &  0.57 & 0.56\\ \hline
0 & Duration &22.51 &  22.50 &   22.50 &  22.51 &  22.51  \\
& Criticality & 0.66 & 0.62 & 0.62 &  0.61 &   0.62\\\hline
0.2 & Duration &22.34& 22.33 &  22.33 &  22.34 &  22.34 \\
& Criticality & 0.64 &   0.63  &0.62 &  0.62 &  0.63\\\hline
0.4 & Duration &
22.09 &   22.10 &  22.10 &  22.10 &  22.10\\
& Criticality &0.62 &  0.60 &  0.60 &  0.61 &  0.61\\\hline
0.6 & Duration &21.46 &  21.48 &  21.48 &  21.48 &  21.48 \\
& Criticality & 0.55 &  0.55 &  0.55 &  0.55 &  0.55 
\\\hline
0.8 & Duration & 20.57 &   20.60 &  20.61 &  20.61 &  20.61 \\
& Criticality & 0.50 &  0.50 &0.50 &  0.48 &   0.50\\\hline
1 & Duration &19.0 &  19.0 &  19.0  &19.0 & 19.0\\
& Criticality &0 & 0&  0& 0& 0 \\
\hline
\end{tabular}
\caption{Triangular activity durations: Worst-case expected project duration and criticality of activity (1,3) with given lower bound on correlation and number  of support points.}  \label{table:triangular}
\end{table}
\subsection{Multi-newsvendor}
\begin{figure}
\begin{center}
\includegraphics[scale=0.6]{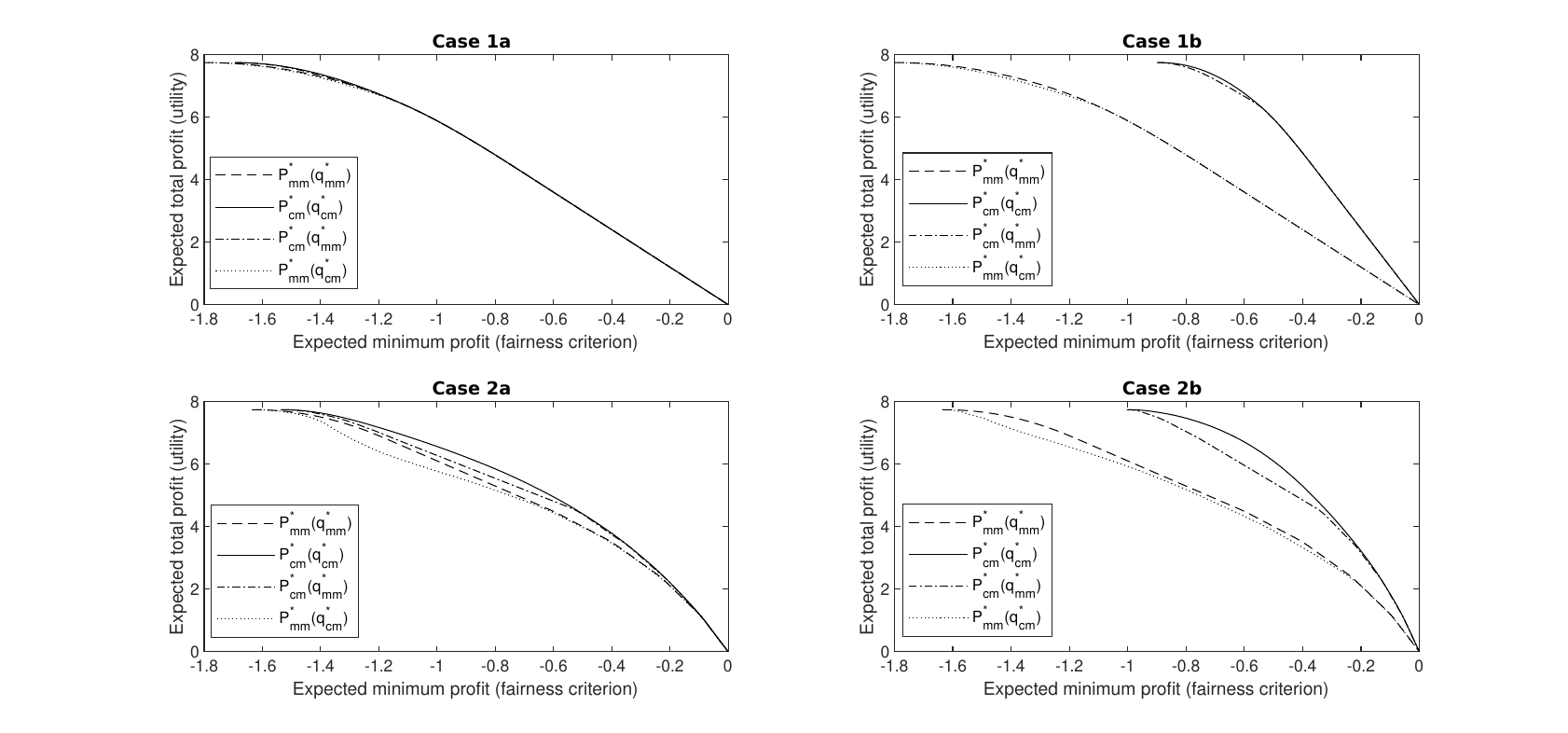}
\caption{Fairness-utility tradeoff for all $\alpha \in [0,1]$ where $\mathbb{P}^*_{\text{mm}}(\mb{q}_{\text{mm}}^*)$ and $\mathbb{P}^*_{\text{cm}}(\mb{q}_{\text{cm}}^*)$ correspond to results with the ambiguity set being correctly specified and $\mathbb{P}^*_{\text{cm}}(\mb{q}_{\text{mm}}^*)$ and $\mathbb{P}^*_{\text{mm}}(\mb{q}_{\text{cm}}^*)$ correspond to results with the ambiguity set being mis-specified.}
\label{fig:bounds}
\end{center}
\end{figure}

\begin{figure}
\begin{center}
\includegraphics[scale=0.5]{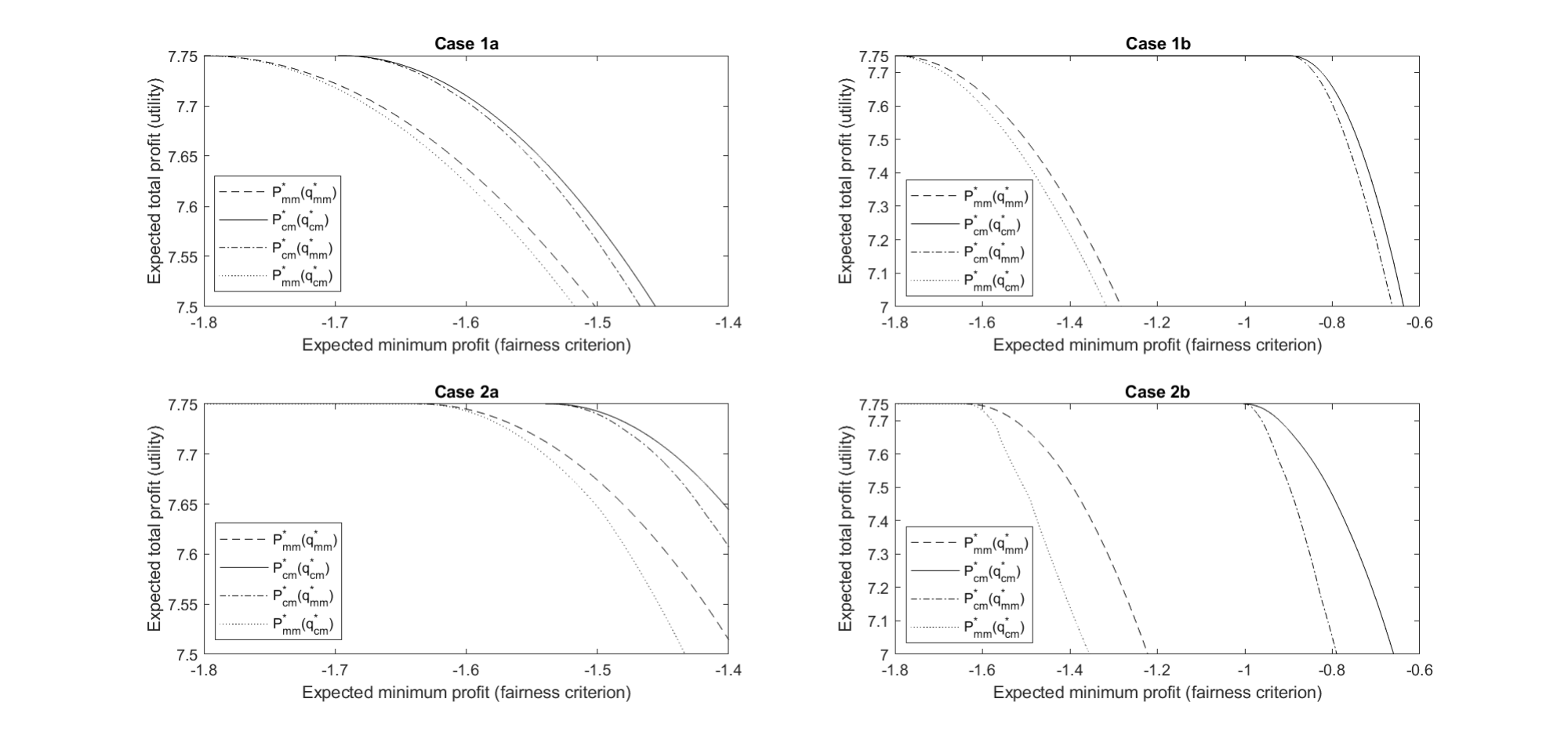}
\caption{Fairness-utility tradeoff for subset of values $\alpha \in [0,1]$ where $\mathbb{P}^*_{\text{mm}}(\mb{q}_{\text{mm}}^*)$ and $\mathbb{P}^*_{\text{cm}}(\mb{q}_{\text{cm}}^*)$ correspond to results with the ambiguity set being correctly specified and $\mathbb{P}^*_{\text{cm}}(\mb{q}_{\text{mm}}^*)$ and $\mathbb{P}^*_{\text{mm}}(\mb{q}_{\text{cm}}^*)$ correspond to results with the ambiguity set being mis-specified.}
\label{fig:bounds00}
\end{center}
\end{figure}

We consider a multi-newsvendor model where a retailer owns $N$ stores. The demand for the product at store $i \in [N]$ is random and is denoted by $\tilde{\xi}_i$. The retailer makes centralized decisions for the stores where the order/production quantity for store $i \in [N]$ is set to $q_i$ and decided before demand is observed. The decision vector $\mb{q} =(q_1,\ldots,q_N)$ is to be chosen from a set ${\cal Q} \subseteq \mathbb{R}^N_+$ which might model constraints on the order quantities such as capacity constraints. The unit selling price at location $i$ is given by $p_i$ and the unit order cost is given by $c_i$. The newsvendor profit at location $i$ for order quantity $q_i$ and demand realization $\tilde{\xi}_i = \xi_i$ is given by $\pi_i(q_i,{\xi}_i) := p_i \min({\xi}_i,q_i)-c_iq_i$. In contrast to the existing literature which primarily focuses on maximizing the expected total profit given by $\mathbb{E}_{\mathbb{P}}[\sum_{i \in [N]}\pi_i(q_i,\tilde{\xi}_i)
]$ where $\tilde{\mb{\xi}} \sim \mathbb{P}$, we also consider the effect of incorporating a Rawlsian fairness criterion into the model (see \cite{hooker}). The Rawlsian fairness criterion focuses on maximizing the expected minimum profit across the locations given by $\mathbb{E}_{\mathbb{P}}[\min_{i \in [N]}\pi_i(q_i,\tilde{\xi}_i)
]$. This encourages improving the performance at the worst-performing location. Such considerations are natural when fairness of allocations is a concern and poor performance at any location can lead to liquidation or the closure of a store. We focus on a distributionally robust model where demands might be correlated and the demand distribution lies in an ambiguity set $\tilde{\mb{\xi}} \sim \mathbb{P} \in {\cal P}$.
The problem is formulated as follows:
\begin{equation} \label{dronews}
\begin{array}{rlll}
\displaystyle \max_{\mbs{q} \in {\cal Q}}  \left(\alpha\inf_{\mathbb{P} \in {\cal P}}\mathbb{E}_{\mathbb{P}}\left[\min_{i \in [N]}\pi_i(q_i,\tilde{\xi}_i)
\right] + (1-\alpha)\inf_{\mathbb{P} \in {\cal P}}\mathbb{E}_{\mathbb{P}}\left[\sum_{i \in [N]}\pi_i(q_i,\tilde{\xi}_i)
\right] \right),
\end{array}
\end{equation}
where $\alpha \in [0,1]$ is a fixed parameter that controls the dual objectives of fairness and total profit. When $\alpha = 0$, the order quantities are chosen to maximize the worst-case total expected profit and when $\alpha = 1$, the order quantities are chosen to maximize the worst-case expected minimum profit across locations. We consider two moment ambiguity sets in our numerical experiments, both which are special cases of the submodular ambiguity set. The first is a marginal moment ambiguity set where the means and variances of the nonnegative demands are fixed:
\begin{equation} \label{eq:bimomabcd}
\begin{array}{lll}
{\cal Q}_{\text{mm}} :=  \{\mathbb{P} \in {\cal P}(\mathbb{R}^N_+) \ | \ \mathbb{E}_{\mathbb{P}}[\tilde{\mb{\xi}}] = \mb{\mu}, \mathbb{E}_{\mathbb{P}}[\mbox{diag}(\tilde{\mb{\xi}}\tilde{\mb{\xi}}')] = \mbox{diag}(\mb{\mu}\mb{\mu}'+\mb{\Sigma}) \big\}.
 \end{array}
 \end{equation} The second ambiguity set is the cross moment ambiguity set ${\cal Q}_{\text{cm}}$ in \eqref{eq:bimom}. We solve the distributionally robust problem \eqref{dronews} by using the dual of the formulation in Theorem \ref{thm:affinemomsdp} by varying the values of $\alpha \in [0,1]$ for each of the ambiguity sets ${\cal P} = {\cal Q}_{\text{mm}}$ and  ${\cal P} = {\cal Q}_{\text{cm}}$ and compute the (worst-case) expected total profit and expected minimium profit. Let $\mb{q}_{\text{mm}}^*$ and $\mb{q}_{\text{cm}}^*$ denote the optimal order quantity vectors respectively for a given $\alpha$. We then use these order quantities to compute the (worst-case) expected total profit and (worst-case) expected minimum profit if the ambiguity set was mis-specified. We recreate mis-specification in the ambiguity set as follows: Suppose the order quantity is $\mb{q}_{\text{mm}}^*$ but the true distribution is chosen from the ambiguity set ${\cal Q}_{\text{cm}}$.  We compute the worst-case expected total profit and worst-case expected minimum profit that arises from this mis-specification. Likewise, we perform the computations where the order quantity is $\mb{q}_{\text{cm}}^*$ but the true distribution is chosen from the ambiguity set ${\cal Q}_{\text{mm}}$. In the numerical experiments, we set $N = 10$, $p_i = 2$ and $c_i = 1$ for all $i \in [N]$. We considered four different settings of the means $\mu_i$, standard deviations $\sigma_i$ and correlation information:
\begin{enumerate}
\item[]{(Case 1a)} Set $\mu_i = 1+i/N$ and $\sigma_i = (1+i/N)/2$ for $i \in [N]$. In this case, demands with higher means have higher standard deviations. In ${\cal Q}_{\text{cm}}$, we set the lower bounds on correlations to $0.4$.
\item[] (Case 1b) With the same mean and standard deviations as case (1a), we set lower bounds on correlations to $0.8$ in ${\cal Q}_{\text{cm}}$.
\item[]{(Case 2a)} Set $\mu_i = 1+i/N$ and $\sigma_i = (1+(N-i+1)/N)/2$ for $i \in [N]$ . In this case, demands with higher means have lower standard deviations. We set lower bounds on correlations to $0.4$ in ${\cal Q}_{\text{cm}}$.
\item[]{(Case 2b)} With the same mean and standard deviations as in case (2a), we set lower bounds on correlations to $0.8$ in ${\cal Q}_{\text{cm}}$.
\end{enumerate}
The results are provided in Figure \ref{fig:bounds} where the fairness-utility tradeoff (worst-case expected minimum profit and worst-case expected total profit) is plotted for the four cases as $\alpha$ is varied over $[0,1]$. A snapshot of the plots is also provided in Figure \ref{fig:bounds00} for restricted ranges of $\alpha$ to more clearly highlight some of the variations. The figures illustrate that by modeling correlation information in the submodular moment ambiguity set, we shift the fairness-utility tradeoff to the right; namely we obtain a better fairness-utility tradeoff by incorporating the correlation information as compared to not incorporating the correlation information. The figure shows that the distributionally robust optimal order quantities are indeed dependent on the correlation structure among the demands and if the model is mis-specified, the performance can drop if the wrong optimal order quantities are used. Lastly, higher the correlations, the more is the improvement provided by the cross moment ambiguity set over the marginal ambiguity set in this application.

\section{Conclusion} \label{sec7}
In this paper, we have introduced a submodular ambiguity set and showcased the modeling of uncertainty using it. With discrete and continuous uncertainty, we have developed linear and semidefinite programs that can be used to incorporate univariate marginal information and bivariate dependence information. The paper shows that the submodular ambiguity set is the natural discrete counterpart of the convex ambiguity set and supplements it for continuous uncertainty, both in modeling and computation. There are several research directions to pursue along this line of work, Firstly, it will be interesting to test the use of decision rules to solve two stage optimization problems under uncertainty with submodular ambiguity sets or to find new decision rules in this context. Secondly, can one identify instances of two stage distributionally robust optimization problems with a submodular ambigutity set that are solvable in polynomial time? Lastly exploring new algorithms for other structured submodular ambiguity sets in applications is a natural direction to consider. We leave this for future research.
\section*{Appendix}
\subsection*{Appendix A}
\noindent \textbf{Proof of Theorem \ref{thm:submod}}\\
Let $\Xi := \prod_{i \in [N]}\Xi_i$. Under assumptions (A1) and (A4), we can reformulate \eqref{droiner} as a finite sized linear program:
\begin{equation*}
\begin{array}{rllll}
\rho^* = \displaystyle \max & \displaystyle \sum_{\mbs{\xi} \in \Xi}\left(\max_{k \in [K]}g_k({\mb{\xi}})\right) p(\mb{\xi})& \\
\mbox{s.t.}
& \displaystyle \sum_{\mbs{\xi} \in \Xi}f_l(\mb{\xi})p(\mb{\xi}) \leq \gamma_l,& \forall l \in [L],\\
 & \displaystyle \sum_{\mbs{\xi} \in \Xi}p(\mb{\xi}) = 1,\\
  & \displaystyle p(\mb{\xi}) \geq 0, & \forall \mb{\xi} \in  \Xi,
\end{array}
\end{equation*}
where the decision variables are $p(\mb{\xi}) = \mathbb{P}(\tilde{\mb{\xi}} = \mb{\xi})$ for $\mb{\xi}  \in \Xi$. We call this the primal linear program. It has a polynomial number of constraints and an exponential number of nonnegative decision variables. The dual linear program is formulated as:
\begin{equation*}\label{dualmdm1}
\begin{array}{rllll}
\rho_{d}^* := \displaystyle \min & \displaystyle y_0 + \sum_{l \in [L]}y_l\gamma_l& \\
\mbox{s.t.}
& \displaystyle y_0 + \sum_{l \in [L]}y_lf_l(\mb{\xi}) \geq \max_{k \in [K]}g_k({\mb{\xi}}),& \forall \mb{\xi} \in \Xi,\\
 & \displaystyle y_l \geq 0, &\forall l \in [L],
\end{array}
\end{equation*}
where the decision variables are $y_0$, $y_l$ for $l \in [L]$ and $\rho_{d}^*$ is the optimal dual objective. The dual linear program has a polynomial number of variables and an exponential number of constraints. The primal linear program is feasible under assumption (A3) and the dual linear program is feasible (set $y_l = 0$ for $l \in [L]$ and $y_0 =  \max_{\mbs{\xi}\in \Xi}\max_{k \in [K]}g_k(\mb{\xi})$).  Strong duality of linear programming guarantees that $\rho^*= \rho_{d}^*$ and the optimal value is finite since the primal and dual linear programs are feasible. The separation problem for the dual linear program is given by:\\
\textit{Given scalars $\bar{y}_0$ and $\bar{y}_l \geq 0$ for $l \in [L]$, decide whether
$$\bar{y}_0 + \sum_{l \in [L]}\bar{y}_lf_l(\mb{\xi}) \geq \max_{k \in [K]}g_k({\mb{\xi}}), \forall \mb{\xi} \in \Xi,$$
and if the answer is no, find a separating hyperplane.}\\
This reduces to checking if
$$\bar{y}_0 + \sum_{l \in [L]}\bar{y}_lf_l(\mb{\xi}) \geq g_k({\mb{\xi}}), \forall \mb{\xi} \in \Xi, \forall k \in [K],$$
which in turn reduces to checking if
$$\bar{y}_0 + \min_{\mbs{\xi} \in \Xi}\left(\sum_{l \in [L]}\bar{y}_lf_l(\mb{\xi})-g_k({\mb{\xi}})\right) \geq 0, \forall k \in [K].$$
Since $\bar{y}_l \geq 0$ for all $l \in [L]$ and $f_l$ is a submodular function for each $l \in [L]$ (assumption (A2)) and $g_k({\mb{\xi}})$ is a supermodular function (assumption (A4)), the function $\sum_{l \in [L]}\bar{y}_lf_l(\mb{\xi})-g_k({\mb{\xi}})$ is submodular for each $k \in [K]$; see Appendix B. Hence to solve the separation problem, we need to solve $K$ submodular function minimization problems where all the functions have efficient evaluation oracles. If we find an index $k^*$ and a vector $\mb{\xi}^*$ such that one of the conditions above is violated, namely:
$$\bar{y}_0 + \sum_{l \in [L]}\bar{y}_lf_l(\mb{\xi}^*)-g_{k^*}({\mb{\xi}^*}) < 0,$$
then a separating hyperplane given by:
$${y}_0 + \sum_{l \in [L]}{y}_lf_l(\mb{\xi}^*)-g_{k^*}({\mb{\xi}^*}) \geq 0.$$
Since all these steps can be done in polynomial time, the dual separation problem is solvable in time polynomial in $N$, $K$, $L$, $\max_{i \in [N]} |\Xi_i|$ and the evaluation time of the oracles. From the equivalence of separation and optimization (see \cite{grotschel1}), the dual linear program is solvable in polynomial time. Hence $\rho^*$ is computable in polynomial time.
\qed
\\

\noindent \textbf{Proof of Theorem \ref{thm:submod1}}\\
Let $\Xi := \prod_{i \in [N]}\Xi_i$. The distributionally robust optimization problem is reformulated as:
\begin{equation*}\label{dualmdm2a}
\begin{array}{rllll}
\displaystyle \inf & \displaystyle y_0 + \sum_{l \in [L]}y_l\gamma_l& \\
\mbox{s.t.}
& \displaystyle y_0 + \sum_{l \in [L]}y_lf_l(\mb{\xi}) \geq \max_{k \in [K]}g_k(\mb{x},{\mb{\xi}}),& \forall \mb{\xi} \in \Xi,\\
 & \displaystyle y_l \geq 0, &\forall l \in [L], \\
  & \displaystyle \mb{x} \in {\cal X}, &
\end{array}
\end{equation*}
where the decision variables are $y_0$, $y_l$ for $l \in [L]$ and $\mb{x}$. The separation problem for this formulation program is given by:\\
\textit{Given scalars $\bar{y}_0$ and $\bar{y}_l \geq 0$ for $l \in [L]$ and vector $\bar{\mb{x}}$, decide whether
$$\bar{\mb{x}} \in {\cal X} \mbox{ and }\bar{y}_0 + \sum_{l \in [L]}\bar{y}_lf_l(\mb{\xi}) \geq \max_{k \in [K]}g_k(\bar{\mb{x}},{\mb{\xi}}), \forall \mb{\xi} \in \Xi,$$
and if the answer is no, find a separating hyperplane.}\\
Verifying the feasibility of $\bar{\mb{x}} \in {\cal X}$ and if it is not feasible, providing a separating hyperplane can be done in polynomial time. The separation problem hence reduces to:\\
\textit{Given scalars $\bar{y}_0$, $\bar{y}_l \geq 0$ for $l \in [L]$ and vector $\bar{\mb{x}} \in {\cal X}$, decide whether
$$\bar{y}_0 + \min_{\mbs{\xi} \in \Xi}\left(\sum_{l \in [L]}\bar{y}_lf_l(\mb{\xi})-g_k(\bar{\mb{x}},{\mb{\xi}})\right) \geq 0, \forall k \in [K].$$}\\
\noindent We need to solve a set of $K$ submodular function minimization problems to solve the separation problem. If all the left hand side values are nonnegative, we are done. Else we find an index $k^* \in [K]$ and a vector $\mb{\xi}^* \in \Xi$ such that:
$$\bar{y}_0 + \sum_{l \in [L]}\bar{y}_lf_l(\mb{\xi}^*)-g_{k^*}(\bar{\mb{x}},\mb{\xi}^*) < 0.$$
Since $g_{k^*}({\mb{x}},\mb{\xi}^*)$ is convex in ${\mb{x}}$ for a fixed $\mb{\xi}^*$, a separating hyperplane is given by:
$$y_0 + \sum_{l \in [L]}y_lf_l(\mb{\xi}^*)- g_{k^*}(\bar{\mb{x}},\mb{\xi}^*)-{\nabla_{\bar{\mbs{x}}} g_{k^*}}(\bar{\mb{x}},\mb{\xi}^*)'({\mb{x}}-\bar{\mb{x}}) \geq 0,$$
where ${\nabla_{\bar{\mbs{x}}} g_{k^*}}(\bar{\mb{x}},\mb{\xi}^*)$ is the subgradient of $g_{k^*} ({\mb{x}},\mb{\xi}^*)$ with respect to $\mb{x}$ evaluated at $\bar{\mb{x}}$. Since all steps can be done in polynomial time, the distributionally robust optimization problem is solvable in polynomial time.
\qed
\\

\noindent \textbf{Proof of Theorem \ref{thm:affine}}\\
Let $\rho_u^*$ be optimal value of the linear program \eqref{primalmdm1}. We prove $\rho^* = \rho_u^*$ in two steps by showing $\rho^* \leq \rho_u^*$ and $\rho^* \leq \rho_u^*$.\\
\noindent \textit{Step (1): $\rho^* \leq \rho_u^*$}\\
To show $\rho_u^*$ is a valid upper bound on $\rho^*$, we start with a probabilistic construction of the formulation \eqref{primalmdm1}. Define for each possible realization $\mb{\xi} \in \prod_{i \in [N]} \Xi_i$, the set of indices that attain the maximum in the objective function as follows:
\begin{center}
$\displaystyle  K(\mb{\xi}) = \arg\max\left\{{\mb{\xi}}'\mb{A}_k{\mb{\xi}}+\mb{b}_k'{\mb{\xi}} + b_{k} \ | \ k \in [K]\right\}$.
\end{center}
For each $\mb{\xi}$, $K(\mb{\xi}) \subseteq [K]$ is a singleton where $|K(\mb{\xi})|= 1$ or contain multiple indices where $|K(\mb{\xi})|> 1$. Given a random vector $\tilde{\mb{\xi}}$ with distribution $\mathbb{P}$, let $k({\mb{\xi}}$) be a measurable selection on the set $K({\mb{\xi}})$ (for example, one can define $k({\mb{\xi}})$ as the smallest index in $K(\mb{\xi})$). Define the decision variables as the probabilities:
\begin{equation*}
\begin{array}{rllll}
\lambda_k  = & \displaystyle \mathbb{P}\left(k(\tilde{\mb{\xi}}) = k\right), & \forall k \in [K], \\
\displaystyle \lambda_{i,k}({\xi}_{i}) = & \displaystyle \mathbb{P}\left(\tilde{\xi}_i = \xi_i, k(\tilde{\mb{\xi}}) = k \right), & \forall {\xi}_i \in \Xi_i, \forall i \in [N], \forall k \in [K],\\
\displaystyle \lambda_{i,j,k}({\xi}_{i},{\xi}_j)  = & \displaystyle \mathbb{P}\left(\tilde{\xi}_i = \xi_i, \tilde{\xi}_j = \xi_j, k(\tilde{\mb{\xi}}) = k\right),  &\forall \xi_i \in \Xi_i, \forall \xi_j \in \Xi_j, \forall (i,j) \in [N]_2, \forall k \in [K].
\end{array}
\end{equation*}
Clearly the variables as defined must satisfy the nonnegativity constraints in \eqref{primalmdm1}. The first six sets of constraints in the formulation are derived from conditions that the variables by definition must satisfy:
\begin{enumerate}
\item Total sum of the probabilities of indices being optimal is one:
\begin{equation*}
\begin{array}{lllll}
\displaystyle \sum_{k \in [K]}\mathbb{P}\left(k(\tilde{\mb{\xi}}) = k\right)  = \displaystyle 1.
\end{array}
\end{equation*}
\item Law of total probability for the index being optimal:
\begin{equation*}
\begin{array}{lllll}
\displaystyle \sum_{\xi_i \in \Xi_i}\mathbb{P}\left(\tilde{\xi}_i = \xi_i, k(\tilde{\mb{\xi}}) = k \right)  = \mathbb{P}\left(k(\tilde{\mb{\xi}}) = k\right).
\end{array}
\end{equation*}
\item Law of total expectation for the univariate functions:
\begin{equation*}
\begin{array}{lllll}
\displaystyle \sum_{k \in [K]}\sum_{\xi_i \in \Xi_i}f_{i,l}({\xi}_i)\mathbb{P}\left(\tilde{\xi}_i = \xi_i, k(\tilde{\mb{\xi}}) = k \right)  = \mathbb{E}\left[f_{i,l}(\tilde{\xi}_i)\right] \in [\underline{\gamma}_{i,l},\overline{\gamma}_{i,l}].
\end{array}
\end{equation*}
\item Consistency of bivariate marginal of $(\tilde{\xi}_i,\tilde{\xi}_j)$ with univariate marginal of $\tilde{\xi}_j$ for the event $k(\tilde{\mb{\xi}}) = k$:
\begin{equation*}
\begin{array}{lllll}
\displaystyle \sum_{{\xi}_i \in \Xi_i}\mathbb{P}\left(\tilde{\xi}_i = \xi_i, \tilde{\xi}_j = \xi_j, k(\tilde{\mb{\xi}}) = k\right) =  \displaystyle \mathbb{P}\left(\tilde{\xi}_j = \xi_j, k(\tilde{\mb{\xi}}) = k\right).
\end{array}
\end{equation*}
\item Consistency of bivariate marginal of $(\tilde{\xi}_i,\tilde{\xi}_j)$ with univariate marginal of $\tilde{\xi}_i$ for the event $k(\tilde{\mb{\xi}}) = k$:
\begin{equation*}
\begin{array}{lllll}
\displaystyle \sum_{{\xi}_j \in \Xi_j}\mathbb{P}\left(\tilde{\xi}_i = \xi_i, \tilde{\xi}_j = \xi_j, k(\tilde{\mb{\xi}}) = k\right) =  \displaystyle \mathbb{P}\left(\tilde{\xi}_i = \xi_i, k(\tilde{\mb{\xi}}) = k\right).
\end{array}
\end{equation*}
\item Law of total expectation for the upper bound on the bivariate functions:
\begin{equation*}
\begin{array}{lllll}
\displaystyle \sum_{k \in [K]}\sum_{{\xi}_i \in \Xi_i}\sum_{{\xi}_j \in \Xi_i}f_{i,j,l}(\xi_i,\xi_j)\mathbb{P}\left(\tilde{\xi}_i = \xi_i, \tilde{\xi}_j = \xi_j,k(\tilde{\mb{\xi}}) = k \right)  = \displaystyle \mathbb{E}\left[f_{i,j,l}(\tilde{\xi}_i,\tilde{\xi}_j)\right] \leq \gamma_{i,j,l}.
\end{array}
\end{equation*}
\end{enumerate}
The objective in \eqref{primalmdm1} is obtained by expressing the expected function value in terms of the decision variables:
\begin{equation*}
\begin{array}{rllll}
\displaystyle \mathbb{E}_{\mathbb{P}}\left[\max_{k \in [K]}\left(\tilde{\mb{\xi}}'\mb{A}_k\tilde{\mb{\xi}}+\mb{b}_k'\tilde{\mb{\xi}} + c_{k}\right)\right]  & = & \displaystyle \sum_{k \in [K]}\mathbb{E}_{\mathbb{P}}\left[\left(\tilde{\mb{\xi}}'\mb{A}_k\tilde{\mb{\xi}}+\mb{b}_k'\tilde{\mb{\xi}} + c_{k}\right)|k(\tilde{\mb{\xi}}) = k\right]\mathbb{P}\left(k(\tilde{\mb{\xi}}) = k\right) \\
& = & \displaystyle \sum_{k \in [K]}\sum_{(i,j) \in [N]_2}\sum_{\xi_i \in \Xi_i}\sum_{\xi_j \in \Xi_j}2{A}_{i,j,k}\xi_i\xi_j\lambda_{i,j,k}({\xi}_{i},{\xi}_j)    \\
& & \displaystyle + \sum_{k \in [K]}\sum_{i \in [N]}\sum_{\xi_i \in {\Xi}_i} \left(A_{i,i,k}\xi_i^2+b_{i,k}\xi_i\right)\lambda_{i,k}(\xi_i)+\sum_{k \in [K]}c_k\lambda_k.
\end{array}
\end{equation*}
From the necessity of all the constraints, we have $\rho^* \leq \rho_u^*$. We next prove sufficiency.\\
\textit{Step (2): $\rho^* \geq \rho_u^*$}\\
We construct a distribution $\mathbb{P}^* \in {\cal P}_{\text{uni}} \cap {\cal P}_{\text{bi-sub}} $ that attains the upper bound $\rho_u^*$ using the optimal solution of the linear program. Consider an optimal solution of the linear program \eqref{primalmdm1} denoted by $\mb{\lambda}^*$. Let ${\cal K}_+ \subseteq [K]$ such that $\lambda_k^* > 0$ for all $k \in {\cal K}_+$ and $\lambda_k ^ * = 0$  for all $k \in {\cal K}_0 := [K]\backslash{\cal K}$. Create a mixture distribution $\mathbb{P}^*$ as follows:\\
(i) Generate a discrete random variable $\tilde{z}$ that takes values in ${\cal K}_+$ with probability $\mathbb{P}^*(\tilde{z} = k) = \lambda_k^*$.\\
(ii) Conditional on the realization of $\tilde{z}$, define the marginal distribution of each random variable ${\tilde{\xi}}_i$ as:
\begin{equation*}
\begin{array}{rlllll}
\displaystyle \mathbb{P}^*\left({\tilde{\xi}_i} = {\xi}_i\big{|} {\tilde{z}} = k\right)  =  \displaystyle \frac{\lambda_{i,k}^{*}(\xi_i)}{\sum_{\xi \in \Xi_i}\lambda_{i,k}^{*}(\xi)},  \forall \xi_i \in \Xi_i.
\end{array}
\end{equation*}
\indent Generate in step (ii), a comonotonic random vector using these conditional marginal distributions. The constraints guarantee that $\sum_{\xi \in \Xi_i}\lambda_{i,k}^{*}(\xi)  > 0$ for all $i \in [N]$, $k \in {\cal K}_+$ and $\lambda_{i,k}^{*}(\xi) = 0$ for all $\xi \in \Xi_i$, $i \in [N]$, $k \in {\cal K}_0$. It is a valid probability measure since $\sum_{k \in {\cal K}_+} \lambda_k^* = 1$. The expected value of $\mathbb{E}[f_{i,l}(\tilde{\xi}_i)]$ in the mixture distribution $\mathbb{P}^*$ is given by:
\begin{equation*}
\begin{array}{rllll}
\displaystyle \mathbb{E}_{\mathbb{P}^*}[f_{i,l}(\tilde{\xi}_i)]  & = & \displaystyle  \sum_{k \in {\cal K}_+} \lambda_k^* f_{i,l}({\xi}_i)\mathbb{P}^*\left({\tilde{\xi}_i} = \xi_i\big{|} {\tilde{z}} = k\right)\\
& = & \displaystyle \sum_{k \in {\cal K}_+} \lambda_k^*\sum_{\xi_i \in \Xi_i}f_{i,l}({\xi}_i) \left(\frac{\lambda_{i,k}^{*}(\xi_i)}{\sum_{\xi \in \Xi_i}\lambda_{i,k}^{*}(\xi)}\right)\\
& = & \displaystyle \sum_{k \in {\cal K}_+} \sum_{\xi_i \in \Xi_i}f_{i,l}({\xi}_i) \lambda_{i,k}^{*}(\xi_i)\\
& &  [\mbox{since } \sum_{\xi \in \Xi_i}\lambda_{i,k}^*(\xi) = \lambda_k^*]\\
& \in & \displaystyle [\underline{\gamma}_{i,l},\overline{\gamma}_{i,l}].
\end{array}
\end{equation*}
Hence the univariate marginal information of $\mathbb{P}^*$ matches the univariate marginal information specified in $ {\cal P}_{\text{uni}}$. Let $\mathbb{Q}_{i,j}^*$ denote a bivariate distribution of the random variables $(\tilde{\xi}_i,\tilde{\xi_j})$ where conditional on $\tilde{z} = k$:
\begin{equation*}
\begin{array}{rlllll}
\displaystyle \mathbb{Q}_{i,j}^*\left({\tilde{\xi}_i} = {\xi}_i,\tilde{\xi}_j = {\xi}_j\big{|} {\tilde{z}} = k\right)  =  \displaystyle \frac{\lambda_{i,j,k}^{*}(\xi_i,\xi_j)}{\sum_{\xi \in \Xi_i}\sum_{\eta \in \Xi_j}\lambda_{i,j,k}^{*}(\xi,\eta)},  \forall \xi_i \in \Xi_i, \forall \xi_j \in \Xi_j, \forall k \in {\cal K}_+.
\end{array}
\end{equation*}
The constraints guarantee that $\sum_{\xi \in \Xi_i}\sum_{\eta \in \Xi_j}\lambda_{i,j,k}^{*}(\xi,\eta)  > 0$ for all $(i,j) \in [N]_2$, $k \in {\cal K}_+$ and $\lambda_{i,j,k}^{*}(\xi,\eta) = 0$ for all $\xi \in \Xi_i$, $\eta \in \Xi_j$, $(i,j) \in [N]_2$, $k \in {\cal K}_0$.
From the feasibility conditions in the linear program, we see that $\mbox{proj}_{i}(\mathbb{Q}_{i,j}^*|k) = \mbox{proj}_{i}(\mathbb{P}^*|k)$ and $\mbox{proj}_{j}(\mathbb{Q}_{i,j}^*|k) = \mbox{proj}_{j}(\mathbb{P}^*|k)$ where $|k$ denotes conditional on $\tilde{z} = k$. This implies the existence of conditional bivariate distributions for each $(\tilde{\xi}_i,\tilde{\xi_j})$ consistent with the conditional marginal distributions of $\tilde{\xi}_i$ and $\tilde{\xi}_j$. However $\mathbb{Q}_{i,j}^*|k$ does not have to be the conditional bivariate distribution for $(\tilde{\xi}_i,\tilde{\xi_j})$ in $\mathbb{P}^*$, namely it is possible that $\mathbb{Q}_{i,j}^*|k \neq \mbox{proj}_{i,j}(\mathbb{P}^*|k)$. The expected value of $\mathbb{E}[f_{i,j,l}(\tilde{\xi}_i,\tilde{\xi}_j)]$ in the mixture distribution $\mathbb{P}^*$ satisfies:
\begin{equation*}
\begin{array}{rllll}
\displaystyle \mathbb{E}_{\mathbb{P}^*}\left[f_{i,j,l}(\tilde{\xi}_i,\tilde{\xi}_j)\right] & = & \displaystyle  \sum_{k \in {\cal K}_+} \lambda_k^* \mathbb{E}_{\mathbb{P}^*}\left[f_{i,j,l}(\tilde{\xi}_i,\tilde{\xi}_j)\big{|} {\tilde{z}} = k\right]\\
& \leq & \displaystyle  \sum_{k \in {\cal K}_+} \lambda_k^* \mathbb{E}_{\mathbb{Q}_{i,j}^*}\left[f_{i,j,l}(\tilde{\xi}_i,\tilde{\xi}_j)\big{|} {\tilde{z}} = k\right] \\
& &  [\mbox{since }f_{i,j,l}(\xi_i,\xi_j) \mbox{ is submodular and }\mathbb{P}^*|k \mbox{ and } \mathbb{Q}_{i,j}^*|k \mbox{ have same marginals}]\\
& = & \displaystyle \sum_{k \in {\cal K}_+} \lambda_k^* \sum_{\xi_i \in \Xi_i}\sum_{\xi_j \in \Xi_j} \left(f_{i,j,l}(\xi_i,\xi_j) \frac{\lambda_{i,j,k}^*(\xi_i,\xi_j)}{\sum_{\xi \in \Xi_i}\sum_{\eta \in \Xi_j}\lambda_{i,j,k}^{*}(\xi,\eta)}\right)\\
& = & \displaystyle \sum_{k \in {\cal K}_+}\sum_{\xi_i \in \Xi_i}\sum_{\xi_j \in \Xi_j} f_{i,j,l}(\xi_i,\xi_j)\lambda_{i,j,k}^*(\xi_i,\xi_j),\\
& & \displaystyle [\mbox{since } \sum_{\xi \in \Xi_i}\sum_{\eta \in \Xi_j}\lambda_{i,j,k}^{*}(\xi,\eta) = \lambda_k^*]\\
& \leq & \displaystyle \gamma_{i,j,l}.
\end{array}
\end{equation*}
Hence $\mathbb{P}^* \in {\cal P}_{\text{bi-sub}}$. The final step is to show the sharpness of the bound under this distribution as follows:
\begin{equation*}
\begin{array}{rllll}
\displaystyle \rho^* & \geq & \displaystyle\mathbb{E}_{\mathbb{P}^*}\left[\max_{k \in [K]}\left(\tilde{\mb{\xi}}'\mb{A}_k\tilde{\mb{\xi}}+\mb{b}_k'\tilde{\mb{\xi}} + c_{k}\right)\right]  \\
& &  [\mbox{since } \mathbb{P}^* \in {\cal P}_{\text{uni}} \cap {\cal P}_{\text{bi-sub}}]\\
& \geq & \displaystyle \sum_{k \in {\cal K}_+} \lambda_k^* \mathbb{E}_{\mathbb{P}^*}\left[\left(\tilde{\mb{\xi}}'\mb{A}_k\tilde{\mb{\xi}}+\mb{b}_k'\tilde{\mb{\xi}} + c_{k}\right) \big{|} \tilde{z} = k\right]\\
& &  [\mbox{evaluating the expected value at the $k$th piece in step (ii) instead of the optimal piece}]\\
& \geq & \displaystyle \sum_{k \in {\cal K}_+} \lambda_k^* \sum_{(i,j) \in [N]_2}\sum_{\xi_i \in \Xi_i}\sum_{\xi_j \in \Xi_j}\left(\frac{2A_{i,j,k}\xi_i\xi_j\lambda_{i,j,k}^{*}(\xi_i,\xi_j) }{\sum_{\xi \in \Xi_i}\sum_{\eta \in \Xi_j}\lambda_{i,j,k}^{*}(\xi,\eta)}\right) \\
&  & \displaystyle + \sum_{k \in {\cal K}_+} \sum_{i \in [N]} \sum_{\xi_i \in \Xi_i}\lambda_k^*  \left(\frac{(A_{i,i,k}\xi_i^2+b_{i,k}\xi_i)\lambda_{i,k}^{*}(\xi_i) }{\sum_{\xi \in \Xi_i}\lambda_{i,k}^{*}(\xi)}\right) +\sum_{k\in {\cal K}_+}c_k\lambda_k^*\\
& & [\mbox{since } A_{i,j,k} \geq 0 \mbox{ and }\mathbb{P}_{|k}^* \mbox{ and } \mathbb{Q}_{i,j|k}^* \mbox{ have same marginals}]\\
& = & \displaystyle \sum_{k \in {\cal K}_+}\sum_{(i,j) \in [N]_2}\sum_{\xi_i \in \Xi_i}\sum_{\xi_j \in \Xi_j}2{A}_{i,j,k}\xi_i\xi_j\lambda_{i,j,k}^*({\xi}_{i},{\xi}_j)  \\
& & \displaystyle + \sum_{k \in {\cal K}_+}\sum_{i \in [N]}\sum_{\xi_i \in \Xi_i} (A_{i,i,k}\xi_i^2+b_{i,k}\xi_i)\lambda_{i,k}^{*}(\xi_i) +\sum_{k \in {\cal K}_+}c_k\lambda_k^*\\
& &  [\mbox{since }  \sum_{\xi \in \Xi_i}\sum_{\eta \in \Xi_j}\lambda_{i,j,k}^{*}(\xi,\eta) = \lambda_k^* \mbox{ and } \sum_{\xi \in \Xi_i}\lambda_{i,k}^*(\xi) = \lambda_k^*]\\
& = & \displaystyle \rho_u^*. \end{array}
\end{equation*}
From steps (1) and (2), $\rho^* = \rho_u^*$.
\qed
\\

\noindent \textbf{Proof of Theorem \ref{thm:submodcon}}\\
Let $\Xi := \prod_{i \in [N]}[0,B_i]$. We can reformulate \eqref{droiner} as a semi-infinite linear program:
\begin{equation*}\label{largeprimalmdm1nops}
\begin{array}{rllll}
\rho^* = \displaystyle \sup & \displaystyle \int_{\mbs{\xi} \in \Xi}\max_{k \in [K]}g_k({\mb{\xi}}) d\theta(\mb{\xi})& \\
\mbox{s.t.}
& \displaystyle \int_{\mbs{\xi} \in \Xi}f_l(\mb{\xi})d\theta(\mb{\xi}) \leq \gamma_l,& \forall l \in [L],\\
 & \displaystyle \int_{\mbs{\xi} \in \Xi}d\theta(\mb{\xi}) = 1,\\
  & \displaystyle d\theta(\mb{\xi}) \geq 0, & \forall \mb{\xi} \in  \Xi,
\end{array}
\end{equation*}
where the decision variable is the  probability measure $\theta$ on the set $\Xi$. This is a semi-infinite linear program with a finite number of constraints and infinite number of nonnegative decision variables. The semi-infinite linear program is feasible by assumption (A3') and since the support of the random vector is bounded, $\rho^*$ is finite. Given $\epsilon > 0$, let $\Xi_{i,\epsilon} = \{0,{\epsilon},2\epsilon,\ldots,\lfloor\frac{B_i}{\epsilon}\rfloor\epsilon\}$ for each $i \in [N]$. Let $\Xi_{\epsilon} := \prod_{i \in [N]}\Xi_{i,\epsilon}$. Then $\Xi_{\epsilon}$ forms an $\epsilon$-grid of $\Xi$. Define $\bar{\Xi}_{\epsilon}=  \Xi_{\epsilon} \cup \{\mb{\xi}_1,\ldots,\mb{\xi}_R\} \subseteq \Xi$ where $\{\mb{\xi}_1,\ldots,\mb{\xi}_R\}$ is the support of the known probability measure in ${\cal P}_{\text{sub}}$. Consider a discretized problem:
\begin{equation*}\label{largeprimalmdm1eps}
\begin{array}{rllll}
{\rho}_{\epsilon}^* := \displaystyle \max & \displaystyle \sum_{\mbs{\xi} \in \bar{\Xi}_{\epsilon}}\max_{k \in [K]}g_k({\mb{\xi}}) p(\mb{\xi})& \\
\mbox{s.t.}
& \displaystyle \sum_{\mbs{\xi} \in \bar{\Xi}_{\epsilon}}f_l(\mb{\xi})p(\mb{\xi}) \leq \gamma_l + M \epsilon,& \forall l \in [L],\\
 & \displaystyle \sum_{\mbs{\xi} \in \bar{\Xi}_{\epsilon}}p(\mb{\xi}) = 1,\\
  & \displaystyle p(\mb{\xi}) \geq 0, & \forall \mb{\xi} \in  \bar{\Xi}_{\epsilon},
\end{array}
\end{equation*}
where $M$ is the Lipschitz constant such that $\max_{i \in [N]}|{\xi}_i -{\chi}_i| \leq z$ implies  $|f_l(\mb{\xi})-f_l(\mb{\chi})| \leq Mz$ for all $l$ and $|g_k(\mb{\xi})-g_k(\mb{\chi})| \leq Mz$ for all $k$. This is a linear program with a finite number of constraints and finite number of decision variables. The number of decision variables is however exponential in $N$. The linear program is feasible since $\sum_{r \in [R]}p_r\delta_{\mbs{\xi}_r}$ is a feasible solution of the linear program. Under the assumptions of submodularity of the functions $f_l$ and supermodularity of the functions $g_k$ with efficient evaluation oracles, computing $\rho_{\epsilon}^*$ is possible in time polynomial in $N$, $K$, $L$, $\max_i B_i$, $R$, $\frac{1}{\epsilon}$ and the evaluation time of the oracles; see Appendix B and the proof of Theorem \ref{thm:submod}. 
Now from Carath\'{e}odory theorem, there exists an optimal probability measure for the semi-infinite linear program that is supported on at most $L + 1$ points in $\Xi$ with optimal value $\rho^*$. Round each support point in this optimal probability measure to the closest point in $\bar{\Xi}_{\epsilon} $, namely within an $\epsilon$ distance, merging points if need be. We thus obtain a modified discrete distribution $\hat{\mathbb{P}}^*_{\epsilon}$ with support contained in $\bar{\Xi}_{\epsilon}$ of cardinality at most $L+1$ that is feasible for the discretized problem. Let $\hat{\rho}^*$ denote the expected objective function value for this new distribution. Then from Lipschitz continuity, $\rho^*-M\epsilon \leq \hat{\rho}^* \leq \rho^* +M\epsilon$. Since $\hat{\rho}^* \leq \rho_{\epsilon}^*$, we get $\rho_{\epsilon}^* \geq  \rho^*-M\epsilon$. The result follows by setting $\epsilon$ appropriately.
\qed
\\

\noindent \textbf{Proof of Theorem \ref{thm:affinemomsdp}}\\
Let $\rho_u^*$ be optimal value of the semidefinite program \eqref{primalmdm1momsdpa}. As before, we prove $\rho^* = \rho_u^*$ in two steps by showing $\rho^* \leq \rho_u^*$ and $\rho_u^* \leq \rho^*$. A key difference here is that we use moments as decision variables rather than probabilities and we need to deal with the issue of the support of the random vector being potentially unbounded.\\ 
\noindent \textit{Step (1): $\rho^* \leq \rho_u^*$}\\
To show that $\rho_u^*$ is a valid upper bound on $\rho^*$, we start by providing a probabilistic interpretation of the formulation. Define the decision variables as:
\begin{equation*}
\begin{array}{rllll}
\displaystyle \lambda_k  = & \displaystyle \mathbb{P}\left(k(\tilde{\mb{\xi}}) = k\right), & \forall k \in [K], \\
\displaystyle \lambda_{i,k}  = & \displaystyle \mathbb{E}_{\mathbb{P}}\left[\tilde{\xi}_i|k(\tilde{\mb{\xi}}) = k\right]\mathbb{P}\left(k(\tilde{\mb{\xi}}) = k\right),  &\forall i \in [N], \forall k \in [K],\\
\displaystyle \lambda_{i,i,k}  = & \displaystyle \mathbb{E}_{\mathbb{P}}\left[\tilde{\xi}_i^2|k(\tilde{\mb{\xi}}) = k\right]\mathbb{P}\left(k(\tilde{\mb{\xi}}) = k\right),  &\forall i \in [N], \forall k \in [K],\\
\displaystyle \lambda_{i,j,k}  = & \displaystyle \mathbb{E}_{\mathbb{P}}\left[\tilde{\xi}_i\tilde{\xi}_j|k(\tilde{\mb{\xi}}) = k\right]\mathbb{P}\left(k(\tilde{\mb{\xi}}) = k\right),  &\forall (i,j) \in [N]_2, \forall k \in [K].
\end{array}
\end{equation*}
The six set of constraints in the formulation are derived from necessary conditions that the variables must satisfy :
\begin{enumerate}
\item Total sum of the probabilities of the indices being optimal is $1$:
\begin{equation*}
\begin{array}{lllll}
\displaystyle \sum_{k \in [K]}\mathbb{P}\left(k(\tilde{\mb{\xi}}) = k\right)  = \displaystyle 1.
\end{array}
\end{equation*}
\item Law of total expectation for the first marginal moment:
\begin{equation*}
\begin{array}{lllll}
\displaystyle \sum_{k \in [K]}\mathbb{E}_{\mathbb{P}}\left[\tilde{\xi}_i|k(\tilde{\mb{\xi}}) = k \right]\mathbb{P}\left(k(\tilde{\mb{\xi}}) = k\right)  = \mathbb{E}_{\mathbb{P}}\left[\tilde{\xi}_i\right] = \mu_i.
\end{array}
\end{equation*}
\item Law of total expectation for the second marginal moment:
\begin{equation*}
\begin{array}{lllll}
\displaystyle \sum_{k \in [K]}\mathbb{E}_{\mathbb{P}}\left[\tilde{\xi}_i^2|k(\tilde{\mb{\xi}}) = k \right]\mathbb{P}\left(k(\tilde{\mb{\xi}}) = k\right)  = \mathbb{E}_{\mathbb{P}}\left[\tilde{\xi}_i^2\right] = \mu_i^2+\Sigma_{i,i}.
\end{array}
\end{equation*}
\item Law of total expectation for the cross moment of $\tilde{\xi}_i$ and $\tilde{\xi}_j$ for $i \neq j$:
\begin{equation*}
\begin{array}{lllll}
\displaystyle \sum_{k \in [K]}\mathbb{E}_{\mathbb{P}}\left[\tilde{\xi}_i\tilde{\xi}_j|k(\tilde{\mb{\xi}}) = k \right]\mathbb{P}\left(k(\tilde{\mb{\xi}}) = k\right)  = \mathbb{E}_{\mathbb{P}}\left[\tilde{\xi}_i\tilde{\xi}_j\right] \geq \mu_i\mu_j+\Sigma_{i,j}.
\end{array}
\end{equation*}
\item Positive semidefiniteness of the three by three moment matrices comes from:
\begin{equation*}
\begin{array}{llll}
\mathbb{E}_{\mathbb{P}}\left[\mathds{1}_{k(\tilde{\mbs{\xi}}) = k}
 \begin{pmatrix}
    1 & \tilde{\xi}_i & \tilde{\xi}_j\\
    \tilde{\xi}_i & \tilde{\xi}_i^2 & \tilde{\xi}_i \tilde{\xi}_j\\
    \tilde{\xi}_j & \tilde{\xi}_i \tilde{\xi}_j &\tilde{\xi}_j^{2}\\
  \end{pmatrix}\right]  \\
   =  \mathbb{P}\left(k(\tilde{\mb{\xi}}) = k\right)\begin{pmatrix}
   1  & \mathbb{E}_{\mathbb{P}}\left[\tilde{\xi}_i|k(\tilde{\mb{\xi}}) = k \right] &\mathbb{E}_{\mathbb{P}}\left[\tilde{\xi}_j|k(\tilde{\mb{\xi}}) = k \right]\\
\mathbb{E}_{\mathbb{P}}\left[\tilde{\xi}_i|k(\tilde{\mb{\xi}}) = k \right] & \mathbb{E}_{\mathbb{P}}\left[\tilde{\xi}_i^2|k(\tilde{\mb{\xi}}) = k \right]& \mathbb{E}_{\mathbb{P}}\left[\tilde{\xi}_i\tilde{\xi}_j|k(\tilde{\mb{\xi}}) = k \right]\\
\mathbb{E}_{\mathbb{P}}\left(\tilde{\xi}_j|k(\tilde{\mb{\xi}}) = k \right) & \mathbb{E}_{\mathbb{P}}\left[\tilde{\xi}_i\tilde{\xi}_j|k(\tilde{\mb{\xi}}) = k \right] &\mathbb{E}_{\mathbb{P}}\left[\tilde{\xi}_j^2|k(\tilde{\mb{\xi}}) = k \right]
  \end{pmatrix}\\
   \succeq 0.
   \end{array}
\end{equation*}
\item Finally, since the support is in $\mathbb{R}^N_+$, we obtain the nonnegativity condition on the conditional marginal moment:
\begin{equation*}
\begin{array}{lllll}
\displaystyle \mathbb{E}_{\mathbb{P}}\left[\tilde{\xi}_i|k(\tilde{\mb{\xi}}) = k \right]\mathbb{P}\left(k(\tilde{\mb{\xi}}) = k\right)  \geq 0.
\end{array}
\end{equation*}
\end{enumerate}
The objective function is obtained by expressing the expected value in terms of the decision variables:
\begin{equation*}
\begin{array}{rllll}
\displaystyle \mathbb{E}_{\mathbb{P}}\left[\max_{k \in [K]}\left(\tilde{\mb{\xi}}'\mb{A}_k\tilde{\mb{\xi}}+\mb{b}_k'\tilde{\mb{\xi}} + c_{k}\right)\right]  & = & \displaystyle \sum_{k \in [K]}\mathbb{E}_{\mathbb{P}}\left[\left(\tilde{\mb{\xi}}'\mb{A}_k\tilde{\mb{\xi}}+\mb{b}_k'\tilde{\mb{\xi}} + c_{k}\right)|k(\tilde{\mb{\xi}}) = k\right]\mathbb{P}\left(k(\tilde{\mb{\xi}}) = k\right) \\
& = & \displaystyle \sum_{k \in [K]}\sum_{(i,j) \in [N]_2} 2A_{i,j,k}\lambda_{i,j,k}\\
& &\displaystyle + \sum_{k \in [K]}\sum_{i \in [N]} \left(A_{i,i,k}\lambda_{i,i,k}+b_{i,k} \lambda_{i,k}\right)+\sum_{k \in [K]}c_k\lambda_k.
\end{array}
\end{equation*}
From the necessity of all the constraints, we have $\rho^* \leq \rho_u^*$. We next prove sufficiency.\\
\textit{Step (2): $\rho^* \geq \rho_u^*$}\\
Consider an optimal solution of the semidefinite program \eqref{primalmdm1momsdpa} which is denoted by $\mb{\lambda}^*$. Partition $[K] = {\cal K}_+ \cup {\cal K}_0$ such that $\lambda_k^* > 0$ for all $k \in {\cal K}_+$ and $\lambda_k ^* = 0$  for all $k \in {\cal K}_0$. It is possible ${\cal K}_0 = \emptyset$ but ${\cal K}_+ \neq \emptyset$ always. Also if $\lambda_{k}^* = 0$, then $\lambda_{i,k}^* = 0$ from the positive semidefiniteness condition of the matrices. We construct a distribution (case (2a)) or a sequence of distributions (case (2b)) that in the limit attains the upper bound $\rho_u^*$ depending on the structure of the optimal solution $\mb{\lambda}^*$. \\
\textit{Case (2a): Values of $\mb{\lambda}^*$ where $\lambda_{i,k}^* = 0 \Rightarrow \lambda_{i,i,k}^* = 0$ for any $k \in {\cal K}$.}\\
In this case $\lambda_{i,j,k}^* = 0$ for all $(i,j) \in [N]_2$, $k \in {\cal K}_0$ from the positive semidefiniteness condition of the matrices. Hence all the optimal decision variables indexed by $k \in {\cal K}_0$ are zero. Let $\epsilon \in (0,1)$.
Create a distribution $\mathbb{P}_{\epsilon}^*$ as follows:\\
(i) Generate a discrete random variable $\tilde{z}$ that takes each value $k \in {\cal K}_+$ with probability $\mathbb{P}_{\epsilon}^*(\tilde{z} = k) = \lambda_k^*$.\\
(ii) Conditional on the realization of $\tilde{z} = k$, define the distribution of the random vector as:  \begin{equation*}
\begin{array}{llllll}
\displaystyle \mathbb{P}_{\epsilon}\left(\tilde{\xi}_i = \frac{\lambda_{i,k}^*}{\lambda_k^*} - \sqrt{\frac{\lambda_{i,i,k}^*}{\lambda_k^*} -\left(\frac{\lambda_{i,k}^*}{\lambda_k^*}\right)^2}\sqrt{\frac{\epsilon}{1-\epsilon}}, \ \forall i \in [N] \ \Bigg{|}  \ \tilde{z} = k\right)  = 1-\epsilon, &\\
\displaystyle \mathbb{P}_{\epsilon}\left(\tilde{\xi}_i =\frac{\lambda_{i,k}^*}{\lambda_k^*} + \sqrt{\frac{\lambda_{i,i,k}^*}{\lambda_k^*} -\left(\frac{\lambda_{i,k}^*}{\lambda_k^*}\right)^2}\sqrt{\frac{1-\epsilon}{\epsilon}}, \ \forall i \in [N] \ \Bigg{|} \ \tilde{z} = k\right)  = \epsilon. &  
\end{array}
\end{equation*}

\noindent Clearly, this is a valid probability measure since $\sum_{k \in {\cal K}_+}\lambda_{k}^*= 1$. Observe that in step (ii), we generate a two point comonotonic random vector with perfect positive dependence. The conditional moments for $k \in {\cal K}_+$ are given by:
\begin{equation*}
\begin{array}{rllll}
\displaystyle 
 \mathbb{E}_{\mathbb{P}_{\epsilon}^*}\left[{\tilde{\xi}_i}\big{|} {\tilde{z}} = k\right] & = & \displaystyle \frac{\lambda_{i,k}^*}{\lambda_k^*}, & \forall i \in [N], \\
\displaystyle 
 \mathbb{E}_{\mathbb{P}_{\epsilon}^*}\left[{\tilde{\xi}_i^2}\big{|} {\tilde{z}} = k\right] & = & \displaystyle \frac{\lambda_{i,i,k}^*}{\lambda_k^*}, & \forall i \in [N], \\
\mathbb{E}_{\mathbb{P}_{\epsilon}^*}\left[{\tilde{\xi}_i\tilde{\xi}_j}\big{|} {\tilde{z}} = k\right] & = & \displaystyle \left(\frac{\lambda_{i,k}^*}{\lambda_k^*}\right)\left(\frac{\lambda_{j,k}^*}{\lambda_k^*}\right) + \sqrt{\frac{\lambda_{i,i,k}^*}{\lambda_k^*} -\left(\frac{\lambda_{i,k}^*}{\lambda_k^*}\right)^2}\sqrt{\frac{\lambda_{j,j,k}^*}{\lambda_k^*} -\left(\frac{\lambda_{j,k}^*}{\lambda_k^*}\right)^2}& \forall (i,j) \in [N]_2.
\end{array}
\end{equation*}
The mean of $\tilde{\xi}_i$ in the distribution $\mathbb{P}_{\epsilon}^*$ is given by:
\begin{equation*}
\begin{array}{rllll}
\displaystyle \mathbb{E}_{\mathbb{P}_{\epsilon}^*}\left[\tilde{\xi}_i\right]  =  \displaystyle  \sum_{k \in {\cal K}_+} \lambda_{k}^* \mathbb{E}_{\mathbb{P}_{\epsilon}^*}\left[{\tilde{\xi}_i}\big{|} {\tilde{z}} = k\right]
 =  \displaystyle \sum_{k \in {\cal K}_+}\lambda_{i,k}^* 
 =  \displaystyle \mu_{i}.
\end{array}
\end{equation*}
The second moment of $\tilde{\xi}_i$ in the distribution $\mathbb{P}_{\epsilon}^*$ is given by:
\begin{equation*}
\begin{array}{rllll}
\displaystyle \mathbb{E}_{\mathbb{P}_{\epsilon}^*}\left[\tilde{\xi}_i^2\right] & = & \displaystyle  \sum_{k \in {\cal K}_+} \lambda_{k}^* \mathbb{E}_{\mathbb{P}_{\epsilon}^*}\left[{\tilde{\xi}_i}^2\big{|} {\tilde{z}} = k\right]
= \displaystyle \sum_{k \in {\cal K}_+} \lambda_{i,i,k}^*
 =  \displaystyle \mu_{i}^2+\Sigma_{i,i}.
\end{array}
\end{equation*}
The cross moment of $\tilde{\xi}_i$ and $\tilde{\xi}_j$ in $\mathbb{P}^*$ for any $\epsilon \in (0,1)$ is given by:
\begin{equation*}
\begin{array}{rllll}
\displaystyle \mathbb{E}_{\mathbb{P}_{\epsilon}^*}\left[\tilde{\xi}_i\tilde{\xi}_j\right] & = & \displaystyle  \sum_{k \in {\cal K}_+} \lambda_{k}
^* \mathbb{E}_{\mathbb{P}_{\epsilon}^*}\left[\tilde{\xi}_i\tilde{\xi}_j\big{|} {\tilde{z}} = k\right]\\
& = & \displaystyle  \sum_{k \in {\cal K}_+}\frac{1}{\lambda_k^*}\left(\lambda_{i,k}^*\lambda_{j,k}^*+\sqrt{\lambda_k^*\lambda_{i,i,k}^*-\lambda_{i,k}^{*2}}\sqrt{\lambda_k^*\lambda_{j,j,k}^*-\lambda_{j,k}^{*2}}\right) \\
& \geq & \displaystyle \sum_{k \in {\cal K}_+}\frac{1}{\lambda_{k}^*} \left(\lambda_{i,k}^*\lambda_{j,k}^*+
\sqrt{\left(\lambda_{k}^*\lambda_{i,j,k}^*-\lambda_{i,k}^*\lambda_{j,k}^*\right)^2}\right)\\
& &  [\mbox{from the positive semidefiniteness of the 3 by 3 matrices}]\\
& = & \displaystyle  \displaystyle \sum_{k \in {\cal K}_+}\frac{1}{\lambda_{k}^*} \left(\lambda_{i,k}^*\lambda_{j,k}^*+
\left|\lambda_{k}^*\lambda_{i,j,k}^*-\lambda_{i,k}^*\lambda_{j,k}^*\right|\right)\\
&\geq  & \displaystyle  \displaystyle   \sum_{k \in {\cal K}_+}\frac{1}{\lambda_{k}^*} \left(\lambda_{i,k}^*\lambda_{j,k}^*+
\lambda_{k}^*\lambda_{i,j,k}^*-\lambda_{i,k}^*\lambda_{j,k}^*\right)\\
&=  & \displaystyle  \displaystyle   \sum_{k \in {\cal K}_+}\lambda_{i,j,k}^* \\
& \geq & \displaystyle \mu_{i}\mu_{j}+\Sigma_{i,j}.
\end{array}
\end{equation*}
The realizations of the random vector are nonnegative for any $\epsilon \in (0,\min_{k \in {\cal K}_+}{\lambda_{i,k}^{*2}}/({\lambda_{k}^*\lambda_{i,i,k}^*}))$ where $0/0 := 1$. The final step is to show the bound is attained by the distribution for any $\epsilon$ in this range. This follows from:
\begin{equation*}
\begin{array}{rllll}
\displaystyle \rho^* & \geq & \displaystyle\mathbb{E}_{\mathbb{P}_{\epsilon}^*}\left[\max_{k \in [K]}\left(\tilde{\mb{\xi}}'\mb{A}_k\tilde{\mb{\xi}}+\mb{b}_k'\tilde{\mb{\xi}} + c_{k}\right)\right]  \\
& &  [\mbox{since } \mathbb{P}_{\epsilon}^* \in {\cal Q}_{\text{cm}} \mbox{ for } \epsilon \in (0,\min_{k \in {\cal K}_+}{\lambda_{i,k}^{*2}}/({\lambda_{k}^*\lambda_{i,i,k}^*}))]\\
& \geq & \displaystyle \sum_{k \in {\cal K}_+} \lambda_k^* \mathbb{E}_{\mathbb{P}_{\epsilon}^*}\left[\left(\tilde{\mb{\xi}}'\mb{A}_k\tilde{\mb{\xi}}+\mb{b}_k'\tilde{\mb{\xi}} + c_{k}\right) \big{|} \tilde{z} = k\right],\\
& &  [\mbox{evaluating the expected value at the $k$th piece in step (ii)}]\\
& \geq & \displaystyle \sum_{k \in {\cal K}_+} \sum_{(i,j) \in [N]_2}2A_{i,j,k}\lambda_{i,j,k}^*+\sum_{k \in {\cal K}_+}\sum_{i \in [N]}\left(A_{i,i,k} \lambda_{i,i,k}^* + b_{i,k}\lambda_{i,k}^*\right)+\sum_{k \in {\cal K}_+}c_k\lambda_k^* \\
& &  [\mbox{since } A_{i,j,k} \geq 0 \mbox{ and } \mathbb{E}_{\mathbb{P}_{\epsilon}^*}[\tilde{\xi}_i\tilde{\xi}_j\big{|} {\tilde{z}} = k] \geq {\lambda_{i,j,k}^*}/{\lambda_k^*}]\\
& = & \displaystyle \rho_u^*. \end{array}
\end{equation*}

\noindent \textit{Case (2b): All other values of $\mb{\lambda}^*$}\\
 Let $\epsilon \in (0,1)$. Here, we create a set of distributions that satisfy the moment constraints and attains the bound in the limit.
 Create a distribution $\mathbb{P}_{\epsilon}^*$ as follows:\\
(i) Generate a discrete random variable $\tilde{z}$ that takes each value $k \in {\cal K}_+$ with probability $\mathbb{P}_{\epsilon}^*(\tilde{z} = k) = \lambda_k^*(1-\epsilon)$ and each value $k \in {\cal K}_0$ with probability $\mathbb{P}_{\epsilon}^*(\tilde{z} = k) = \epsilon/|{\cal K}_0|$.\\
(ii) Conditional on the realization of $\tilde{z} = k$, define the distribution of the random vector as follows:\\
If $k \in {\cal K}_+$, then:
\begin{equation*}
\begin{array}{llllll}
\displaystyle \mathbb{P}_{\epsilon}\left(\tilde{\xi}_i = \frac{\lambda_{i,k}^*}{\lambda_k^*}, \ \forall i \in [N] \ \Bigg{|}  \ \tilde{z} = k\right)  = 1-\epsilon, &\\
\displaystyle \mathbb{P}_{\epsilon}\left(\tilde{\xi}_i = \frac{\lambda_{i,k}^*}{\lambda_k^*}+ \frac{1}{\sqrt{\epsilon}}\sqrt{{\frac{\lambda_{i,i,k}^*}{\lambda_k^*}-\left(\frac{\lambda_{i,k}^*}{\lambda_k^*}\right)^2}}, \ \forall i \in [N] \ \Bigg{|} \ \tilde{z} = k\right)  = \epsilon. &  
\end{array}
\end{equation*}
If $k \in {\cal K}_0$, then:
\begin{equation*}
\begin{array}{llllll}
\displaystyle \mathbb{P}_{\epsilon}\left(\tilde{\xi}_i = \sqrt{\frac{\lambda_{i,i,k}^*|{\cal K}_0|}{\epsilon}}, \ \forall i \in [N] \ \Bigg{|}  \ \tilde{z} = k\right)  = 1.
\end{array}
\end{equation*}

\noindent Clearly, it is a valid probability measure since $\sum_{k \in {\cal K}_+ \cup {\cal K}_0}\lambda_{k}^*= 1$.  Observe that in step (ii), we generate a two point comonotonic random vector when $\tilde{z} = k \in {\cal K}_+$ and a one point random vector when $\tilde{z} = k \in {\cal K}_0$. It is straightforward to verify that the conditional moments in the distribution are given for each $k \in {\cal K}_+$ by:
\begin{equation*}
\begin{array}{rllll}
\displaystyle 
 \mathbb{E}_{\mathbb{P}_{\epsilon}^*}\left[{\tilde{\xi}_i}\big{|} {\tilde{z}} = k\right] & = & \displaystyle \frac{\lambda_{i,k}^*}{\lambda_k^*} +\sqrt{\epsilon}\sqrt{\frac{\lambda_{i,i,k}^*}{\lambda_k^*} -\left(\frac{\lambda_{i,k}^*}{\lambda_k^*}\right)^2}, & \forall i \in [N], \\
\displaystyle 
 \mathbb{E}_{\mathbb{P}_{\epsilon}^*}\left[{\tilde{\xi}_i^2}\big{|} {\tilde{z}} = k\right] & = & \displaystyle \frac{\lambda_{i,i,k}^*}{\lambda_k^*}+2\sqrt{\epsilon}\left(\frac{\lambda_{i,k}^*}{\lambda_k^*}\right)\sqrt{\frac{\lambda_{i,i,k}^*}{\lambda_k^*} -\left(\frac{\lambda_{i,k}^*}{\lambda_k^*}\right)^2}, & \forall i \in [N], \\
\mathbb{E}_{\mathbb{P}_{\epsilon}^*}\left[{\tilde{\xi}_i\tilde{\xi}_j}\big{|} {\tilde{z}} = k\right] & = & \displaystyle \left(\frac{\lambda_{i,k}^*}{\lambda_k^*}\right)\left(\frac{\lambda_{j,k}^*}{\lambda_k^*}\right) + \sqrt{\frac{\lambda_{i,i,k}^*}{\lambda_k^*} -\left(\frac{\lambda_{i,k}^*}{\lambda_k^*}\right)^2}\sqrt{\frac{\lambda_{j,j,k}^*}{\lambda_k^*} -\left(\frac{\lambda_{j,k}^*}{\lambda_k^*}\right)^2}&  \\
&& \displaystyle +\sqrt{\epsilon}\left(\left(\frac{\lambda_{i,k}^*}{\lambda_k^*}\right)\sqrt{\frac{\lambda_{j,j,k}^*}{\lambda_k^*} -\left(\frac{\lambda_{j,k}^*}{\lambda_k^*}\right)^2}+\left(\frac{\lambda_{j,k}^*}{\lambda_k^*}\right)\sqrt{\frac{\lambda_{i,i,k}^*}{\lambda_k^*} -\left(\frac{\lambda_{i,k}^*}{\lambda_k^*}\right)^2}\right), & \forall (i,j) \in [N]_2,\\
\end{array}
\end{equation*}
and for each $k \in {\cal K}_0$ by:
\begin{equation*}
\begin{array}{rllll}
\displaystyle 
 \mathbb{E}_{\mathbb{P}_{\epsilon}^*}\left[{\tilde{\xi}_i}\big{|} {\tilde{z}} = k\right] & = & \displaystyle \sqrt{\frac{\lambda_{i,i,k}^*|{\cal K}_0|}{\epsilon}}, & \forall i \in [N], \\
\displaystyle 
 \mathbb{E}_{\mathbb{P}_{\epsilon}^*}\left[{\tilde{\xi}_i^2}\big{|} {\tilde{z}} = k\right] & = & \displaystyle {\frac{\lambda_{i,i,k}^*|{\cal K}_0|}{\epsilon}}, & \forall i \in [N], \\
\mathbb{E}_{\mathbb{P}_{\epsilon}^*}\left[{\tilde{\xi}_i\tilde{\xi}_j}\big{|} {\tilde{z}} = k\right] & = & \displaystyle \sqrt{\lambda_{i,i,k}^*\lambda_{j,j,k}^*}\frac{|{\cal K}_0|}{\epsilon}& \forall (i,j) \in [N]_2.
\end{array}
\end{equation*}
As $\epsilon \downarrow 0$, the mean of $\tilde{\xi}_i$ converges to:
\begin{equation*}
\begin{array}{rllll}
\displaystyle \lim_{\epsilon \downarrow 0}\mathbb{E}_{\mathbb{P}_{\epsilon}^*}\left[\tilde{\xi}_i\right] & = & \displaystyle  \lim_{\epsilon \downarrow 0}\left(\sum_{k \in {\cal K}_+} \lambda_{k}^*(1-\epsilon) \mathbb{E}_{\mathbb{P}_{\epsilon}^*}\left[{\tilde{\xi}_i}\big{|} {\tilde{z}} = k\right]+\sum_{k \in {\cal K}_0} \frac{\epsilon}{|K_0|}\mathbb{E}_{\mathbb{P}_{\epsilon}^*}\left[{\tilde{\xi}_i}\big{|} {\tilde{z}} = k\right]\right)\\
& = & \displaystyle \sum_{k \in {\cal K}_+}\lambda_{i,k}^*\\
& =& \displaystyle \mu_i. 
\end{array}
\end{equation*}
As $\epsilon \downarrow 0$, the second moment of $\tilde{\xi}_i$ converges to:
\begin{equation*}
\begin{array}{rllll}
\displaystyle \lim_{\epsilon \downarrow 0}\mathbb{E}_{\mathbb{P}_{\epsilon}^*}\left[\tilde{\xi}_i^2\right] & = &   \displaystyle  \lim_{\epsilon \downarrow 0}\left(\sum_{k \in {\cal K}_+} \lambda_{k}^*(1-\epsilon) \mathbb{E}_{\mathbb{P}_{\epsilon}^*}\left[{\tilde{\xi}_i^2}\big{|} {\tilde{z}} = k\right]+\sum_{k \in {\cal K}_0} \frac{\epsilon}{|K_0|}\mathbb{E}_{\mathbb{P}_{\epsilon}^*}\left[{\tilde{\xi}_i^2}\big{|} {\tilde{z}} = k\right]\right)\\
& = & \displaystyle \sum_{k \in {\cal K}_+} \lambda_{i,i,k}^* +\sum_{k \in {\cal K}_0}\lambda_{i,i,k}^{*}\\
& = & 
 \displaystyle \mu_{i}^2+\Sigma_{i,i}.
\end{array}
\end{equation*}
As $\epsilon \downarrow 0$, the cross moment of $\tilde{\xi}_i$ and $\tilde{\xi}_j$ satisfies the condition:
\begin{equation*}
\begin{array}{rllll}
\displaystyle \lim_{\epsilon \downarrow 0}\mathbb{E}_{\mathbb{P}_{\epsilon}^*}\left[\tilde{\xi}_i\tilde{\xi}_j\right] & = & \displaystyle  \lim_{\epsilon \downarrow 0}\left(\sum_{k \in {\cal K}_+} \lambda_{k}
^*(1-\epsilon) \mathbb{E}_{\mathbb{P}_{\epsilon}^*}\left[\tilde{\xi}_i\tilde{\xi}_j\big{|} {\tilde{z}} = k\right] + \sum_{k \in {\cal K}_0}\frac{\epsilon}{|{\cal K}_0|}\mathbb{E}_{\mathbb{P}_{\epsilon}^*}\left[\tilde{\xi}_i\tilde{\xi}_j\big{|} {\tilde{z}} = k\right]\right)\\
& = & \displaystyle  \sum_{k \in {\cal K}_+}\frac{1}{\lambda_k^*}\left(\lambda_{i,k}^*\lambda_{j,k}^*+\sqrt{\lambda_k^*\lambda_{i,i,k}^*-\lambda_{j,k}^{*2}}\sqrt{\lambda_k^*\lambda_{j,j,k}^*-\lambda_{j,k}^{*2}}\right)  + \sum_{k \in {\cal K}_0}\sqrt{\lambda_{i,i,k}^*\lambda_{j,j,k}^*}\\\
& \geq & \displaystyle \sum_{k \in {\cal K}_+}\frac{1}{\lambda_{k}^*} \left(\lambda_{i,k}^*\lambda_{j,k}^*+
\sqrt{\left(\lambda_{k}^*\lambda_{i,j,k}^*-\lambda_{i,k}^*\lambda_{j,k}^*\right)^2}\right)+ \sum_{k \in {\cal K}_0} \sqrt{\lambda_{i,j,k}^{*2}}\\
& & [\mbox{from the positive semidefiniteness of the 3 by 3 matrices}]\\
& \geq  & \displaystyle  \displaystyle   \sum_{k \in {\cal K}}\lambda_{i,j,k}^* \\
& \geq & \displaystyle \mu_{i}\mu_{j}+\Sigma_{i,j}.
\end{array}
\end{equation*}
The support of the random vector is nonnegative by construction. The final step is to show the sharpness of the bound. As $\epsilon \downarrow 0$, we get:
\begin{equation*}
\begin{array}{rllll}
\displaystyle \rho^* & \geq & \displaystyle \lim_{\epsilon \downarrow 0}\mathbb{E}_{\mathbb{P}_{\epsilon}^*}\left[\max_{k \in [K]}\left(\tilde{\mb{\xi}}'\mb{A}_k\tilde{\mb{\xi}}+\mb{b}_k'\tilde{\mb{\xi}} + c_{k}\right)\right]  \\
& &  [\mbox{using the limit}]\\
& \geq & \displaystyle \lim_{\epsilon \downarrow 0}\sum_{k \in {\cal K}_+} \lambda_k^*(1-\epsilon) \mathbb{E}_{\mathbb{P}_{\epsilon}^*}\left[\left(\tilde{\mb{\xi}}'\mb{A}_k\tilde{\mb{\xi}}+\mb{b}_k'\tilde{\mb{\xi}} + c_{k}\right) \big{|} \tilde{z} = k\right] + \lim_{\epsilon \downarrow 0}\sum_{k \in {\cal K}_0}\frac{\epsilon}{|{\cal K}_0|}\mathbb{E}_{\mathbb{P}_{\epsilon}^*}\left[\left(\tilde{\mb{\xi}}'\mb{A}_k\tilde{\mb{\xi}}+\mb{b}_k'\tilde{\mb{\xi}} + c_{k}\right) \big{|} \tilde{z} = k\right],\\
& &  [\mbox{evaluating the expected value at the $k$th piece in step (ii)}]\\
& \geq & \displaystyle \sum_{k \in {\cal K}} \sum_{(i,j) \in [N]_2}2A_{i,j,k}\lambda_{i,j,k}^*+\sum_{k \in {\cal K}}\sum_{i \in [N]}\left(A_{i,i,k} \lambda_{i,i,k}^* + b_{i,k}\lambda_{i,k}^*\right)+\sum_{k \in {\cal K}}c_k\lambda_k^* \\
& &  [\mbox{since } A_{i,j,k} \geq 0]\\
& = & \displaystyle \rho_u^*. \end{array}
\end{equation*}
From steps (1) and (2), $\rho^* = \rho_u^*$.
\qed
\subsection*{Appendix B: Review of submodularity}
We list out three examples of submodular functions that satisfy the definition in \eqref{submodular1}:
\texitem{(i)} $f(\mb{\xi}) = h(\mb{a}'\mb{\xi})$ where $\mb{a} \geq \mb{0}$ and $h: \mathbb{R} \rightarrow \mathbb{R}$ is a concave function,
\texitem{(ii)} $f(\mb{\xi}) = \max(\xi_1,\ldots,\xi_N)$,
\texitem{(iii)} $f(\mb{\xi}) = -\prod_{i \in [N]} \xi_i$ where $\mb{\xi} \geq \mb{0}$.

\noindent As these examples show, submodular functions might be convex functions (example (ii)), concave functions (example (i)), neither convex or concave (example (iii)) or both convex and concave (the linear function $f(\mb{\xi}) = \mb{a}'\mb{\xi}$). Other examples of submodular functions include the weighted cut function in a directed graph, the entropy of a random vector, the influence function in a social network and the rank function of a matroid (see \cite{Lovasz,topkis,bach} for examples of submodular and supermodular functions arising in graph theory, probability, operations research, game theory, machine learning and artificial intelligence). While submodularity is known to be preserved under certain operations such as taking nonnegative weighted sum of submodular functions or taking a partial minimum of a submodular function, it is not preserved under operations such as taking the pointwise maximum or minimum of submodular functions. Univariate functions are both submodular and supermodular and hence all (additively) separable functions are both submodular and supermodular (see Figure \ref{fig:sub} for an illustration). Bivariate functions which are convex functions of the difference of the two variables is also submodular.

\begin{figure}
\begin{center}
  \begin{tikzpicture}[scale=0.8]
\tikzset{vertex/.style = {shape=circle,draw,minimum size=0.2em}}
\tikzset{edge/.style = {->,> = latex}}
 \node (1) at (-0.5,1) {};
  \fill (1) circle (0) node [left] {\mbox{Submodular}};
  \draw[rotate around={-45:(0,0)}] (0,0) ellipse (130 pt and 60pt);
     \node (3) at (8,1) {};
  \fill (3) circle (0) node [rotate=0,left] {\mbox{Supermodular}};
  \draw[rotate around={45:(4.5,0)}] (4.5,0) ellipse (130 pt and 60pt);
   \node (2) at (-1.5,4) {};
  \fill (2) circle (0) node [right] {\mbox{Convex}};
  \draw[rotate around={-45:(1.5,2)}] (1.5,2) ellipse (130 pt and 60pt);
   \node (4) at (3,4) {};
  \fill (4) circle (0) node [rotate=0,right] {\mbox{Concave}};
  \draw[rotate around={45:(2.8,2)}] (2.8,2) ellipse (130 pt and 60pt);
     \node (5) at (1,-2.5) {};
  \fill (5) circle (0) node [rotate=0,right] {\mbox{Separable}};
     \node (6) at (1.4,-0.5) {};
  \fill (6) circle (0) node [rotate=0,right] {\mbox{Affine}};
   \node (7) at (-0.3,-1) {};
    \fill (7) circle (0.1) node [above] {(i)};
     \node (8) at (-0.3,2) {};
    \fill (8) circle (0.1) node [above] {(ii)};
       \node (9) at (-2,0) {};
    \fill (9) circle (0.1) node [above] {(iii)};
\end{tikzpicture}
\caption{Venn diagram illustrating the set of convex, concave, submodular and supermodular functions defined over $\mathbb{R}^N_+$. Affine functions lie at the intersection of all four sets.}
\label{fig:sub}
\end{center}
\end{figure}
Submodular functions defined over discrete domains behave somewhat similarly to convex functions defined over continuous domains, particularly in terms of the minimization of such functions (see \cite{Lovasz}). Consider the submodular function minimization problem:
\begin{equation} \label{SFM}
\begin{array}{rlll}
\displaystyle \min_{\mbs{\xi} \in \prod_{i \in [N]}\Xi_i} f(\mb{\xi}).
\end{array}
\end{equation}
Assume the submodular function is given by an polynomial time evaluation oracle, namely given a $\mb{\xi} \in \prod_{i \in [N]} \Xi_i$, the oracle returns $f(\mb{\xi})$ in polynomial time\footnote{When the function is rational valued, the oracle returns it exactly in polynomial time. For real valued functions, the oracle is assumed to return a rational value with relative error $\epsilon > 0$ in time polynomial in the input size and $\log(1/\epsilon$).}. A key result proved in \cite{grotschel1}, building on the work of \cite{edmonds,khachiyan}, showed that \eqref{SFM} is solvable in time polynomial in the input size using the ellipsoid method when the sets $\Xi_i$ are discrete and finite. A key object that aids in the efficient minimization of submodular functions is the Lov\'{a}sz extension \cite{Lovasz} (also known as the Choquet integral \cite{choquet}) which is defined through the construction of a comonotonic random vector. Let ${\cal P}_F(\mathbb{P}_1,\ldots,\mathbb{P}_N)$ denote the Fr\'{e}chet set of distributions where $\tilde{\xi}_i \sim \mathbb{P}_i$ with $\text{supp}(\mathbb{P}_i) \subseteq \Xi_i$ for each $i \in [N]$. Let $F_i(\cdot)$ denote the cumulative distribution function of $\tilde{\xi}_i$ for each $i \in [N]$. A comonotonic random vector has maximal positive dependence in the Fr\'{e}chet set of distributions and is given by:
$$\tilde{\mb{\xi}}^{c} := (F_1^{-1}(\tilde{U}),\ldots,F_N^{-1}(\tilde{U})),$$
where $\tilde{U}$ is a uniform random variable on $[0,1]$ and $F_i^{-1}(\cdot)$ is the generalized inverse distribution function of the cumulative distribution function $F_i(\cdot)$. Let $\mathbb{P}^c$ denote the distribution of the comonotonic random vector.  Clearly $\mathbb{P}^c \in {\cal P}_F(\mathbb{P}_1,\ldots,\mathbb{P}_N)$. The support of the comonotonic random vector is contained in completely ordered subsets of $\mathbb{R}^N$ with:
\begin{equation} \label{comonup}
\begin{array}{rlll}
\displaystyle \mathbb{P}^c(\tilde{\mb{\xi}}^c > \mb{t}) = \min_{i \in [N]}\mathbb{P}_i(\tilde{{\xi}}_i > t_i) \mbox{ and } \mathbb{P}^c(\tilde{\mb{\xi}}^c \leq \mb{t}) = \min_{i \in [N]}\mathbb{P}_i(\tilde{{\xi}}_i \leq {t}_i), \forall \mb{t} \in \mathbb{R}^N.
\end{array}
\end{equation}
The expected value of a function of the comonotonic random vector is computed as:
 \begin{equation} \label{exp}
\begin{array}{rlll}
\displaystyle \mathbb{E}_{\mathbb{P}^c}\left[f(\tilde{\mb{\xi}}^c)\right] = \int_{0}^{1}f\left(F_1^{-1}({t}),\ldots,F_N^{-1}({t})\right)dt.
\end{array}
\end{equation}
  With discrete and finite $\Xi_i$, the cardinality of the support of the comonotonic random vector is at most $\sum_{i \in [N]}|{\Xi}_i|$ and the value in \eqref{exp} is computable using a polynomial number of calls to the evaluation oracle. A well known extremal characterization of the comonotonic random vector (see \cite{tchen,ruschbook}) for general marginals $\mathbb{P}_1,\ldots,\mathbb{P}_N$ is given by:
\begin{equation} \label{droinfer}
\begin{array}{rlll}
\displaystyle \inf_{\mathbb{P} \in {\cal P}_F(\mathbb{P}_1,\ldots,\mathbb{P}_N)}\mathbb{E}_{\mathbb{P}}\left[f(\tilde{\mb{\xi}})\right] = \mathbb{E}_{\mathbb{P}^c}\left[f(\tilde{\mb{\xi}}^c)\right], \forall \mbox{submodular } f \mbox{ such that expectations exists}.
\end{array}
\end{equation}
The extremal characterization in \eqref{droinfer} has been proved to be very useful in developing efficient algorithms for submodular function minimization. Specifically, the function $f: \prod_{i \in [N]}\Xi_i \rightarrow \mathbb{R}$ is submodular if and only if the functional $\inf_{\mathbb{P} \in {\cal P}_F(\mathbb{P}_1,\ldots,\mathbb{P}_N)}\mathbb{E}_{\mathbb{P}}[f(\tilde{\mb{\xi}})]: \mathbb{P}_1,\ldots,\mathbb{P}_N\rightarrow \mathbb{R}$ is convex (see \cite{Lovasz,bach}). The domain of the functional is specified by the marginal probability measures $\mathbb{P}_1,\ldots,\mathbb{P}_N$ where $\mbox{supp}(\mathbb{P}_1) \subseteq \Xi_1,\ldots,\mbox{supp}(\mathbb{P}_N) \subseteq \Xi_N$. The corresponding functional value is exactly the Lov\'{a}sz extension or the Choquet integral. This lets one transform the submodular function minimization problem in \eqref{SFM} to a convex minimization problem over probability measures as follows:
\begin{equation} \label{SFM1}
\begin{array}{rlll}
\displaystyle \inf_{\mbs{\xi} \in \prod_{i \in [N]}\Xi_i} f(\mb{\xi}) \overset{(a)}{=} \displaystyle \inf_{\small{\text{supp}(\mathbb{P}_i)} \subseteq \Xi_i, \forall i \in [N]} \inf_{\mathbb{P} \in {\cal P}_F(\mathbb{P}_1,\ldots,\mathbb{P}_N)} \mathbb{E}_{\mathbb{P}}\left[f(\tilde{\mb{\xi}})\right] \overset{(b)}{=} \inf_{\small{\text{supp}(\mathbb{P}_i)} \subseteq \Xi_i, \forall i \in [N]} \mathbb{E}_{\mathbb{P}^c}\left[f(\tilde{\mb{\xi}}^c)\right].
\end{array}
\end{equation}
The first equality holds since the optimal solution on the right hand side of (a) will be attained at a Dirac measure where the functional value and the function value coincide. Equality (b) comes from \eqref{droinfer}. When the domain is discrete and finite, the convex minimization problem on the right hand side of (b) is solvable in polynomial time using the ellipsoid method where both the function value and its subgradient are computable in polynomial time (see \cite{Lovasz,bach}). .

\renewcommand{\baselinestretch}{1.00}

\end{document}